\documentclass[10pt]{article}
\usepackage{amscd,amsfonts,amssymb,amscd,amsmath,latexsym}
\usepackage[all]{xy}
\makeatletter
\renewcommand{\subsubsection}{\@startsection
{subsubsection}
{1}
{0mm}
{0mm}
{0mm}
{\normalfont\normalsize\itshape}}
\makeatother
  
\newtheorem{theorem}{Theorem}[section] 
\newtheorem{prop}[theorem]{Proposition}
\newtheorem{lem}[theorem]{Lemma}
\newtheorem{ddd}[theorem]{Definition}

\newcommand{\colim}{\mathrm{colim}}

\newcommand{\K}{\mathbb{K}}

\newcommand{\forget}[1]{}

{\catcode`@=11\global\let\c@equation=\c@theorem}

\newcommand{\Z}{\mathbb{Z}}

\newcommand{\diag}{{\tt diag}}
\newcommand{\proof}{{\it Proof.$\:\:\:\:$}}

\newcommand{\R}{\mathbb{R}}

\newcommand{\Q}{\mathbb{Q}}

\newcommand{\bS}{\mathbf{S}}

\newcommand{\bA}{\mathbf{A}}

\newcommand{\ch}{\mathbf{ch}}

\newcommand{\C}{\mathbb{C}}

\newcommand{\Aut}{{\tt Aut}}

\newcommand{\Tr}{{\tt Tr}}
\newcommand{\cT}{\mathcal{T}}

\newcommand{\cZ}{\mathcal{Z}}

\newcommand{\cE}{\mathcal{E}}

\newcommand{\cW}{\mathcal{W}}

\newcommand{\cC}{\mathcal{C}}

\newcommand{\cA}{\mathcal{A}}
\newcommand{\cU}{\mathcal{U}}
\newcommand{\Hom}{{\tt Hom}}

\newcommand{\cO}{\mathcal{O}}
\newcommand{\End}{{\tt End}}
\newcommand{\Ext}{{\tt Ext}}

\newcommand{\cD}{\mathcal{D}}

\newcommand{\tr}{{\tt tr}}

\newcommand{\id}{{\tt id}}

\newcommand{\nat}{\mathbb{N}}
\newcommand{\supp}{{\tt supp}}

\newcommand{\cB}{\mathcal{B}}
\def\imath{{i}}

\def\hB{\hspace*{\fill}$\Box$ \newline\noindent}

\newcommand{\cS}{\mathcal{S}}

\def\hB{\hspace*{\fill}$\Box$ \\[0cm]\noindent}
\newcommand{\cL}{\mathcal{L}}
 \newcommand{\cG}{\mathcal{G}}
 \newcommand{\cX}{\mathcal{X}}

\newcommand{\pr}{{\tt pr}}
\newcommand{\bX}{\mathbf{X}}
\newcommand{\bY}{\mathbf{Y}}

\newcommand{\ev}{{\tt ev}}

\newcommand{\Sets}{\mathcal{S}ets}

\title{Inertia and delocalized twisted cohomology}
\author{Ulrich Bunke, Thomas Schick and Markus Spitzweck
\thanks{Mathematisches Institut, Universit{\"a}t G{\"o}ttingen,
Bunsenstr. 3-5, 37073 G{\"o}ttingen, GERMANY,  bunke@uni-math.gwdg.de,
schick@uni-math.gwdg.de, spitz@uni-mah.gwdg.de}
}

\begin{document}
\newcommand{\Top}{{\tt Top}}
\newcommand{\bM}{{\mathbf{M}}}
\newcommand{\ori}{{\tt or}} 
\newcommand{\cat}{{\tt Cat}}
\newcommand{\cov}{{\tt cov}}
\newcommand{\Ab}{{\tt Ab}}
\newcommand{\Sh}{{\tt Sh}}
\newcommand{\Site}{{\tt Site}}
\newcommand{\Mf}{{\tt Mf}}
\renewcommand{\Sets}{{\tt Sets}}
\renewcommand{\Pr}{{\tt PSh}}
\newcommand{\eins}{{\mathbf{1}}}
\newcommand{\Tw}{{\tt Tw}}
\newcommand{\Tors}{{\tt Tors}}

%\begin{Abstract}
%\end{Abstract}

\maketitle 
\tableofcontents

\newcommand{\bG}{\mathbf{G}}

\section{Introduction}

\subsection{Motivation}

\subsubsection{}

In the recent mathematical literature cohomological and topological properties of orbifolds became an intensively studied subject. A considerable part of the motivation comes from the mirror symmetry program where orbifolds arise naturally. Cornerstones\footnote{Here we mention those works which are relevant for the present paper. Note that there is a huge literature on orbifolds in algebraic geometry and mathematical physics.}  of  the recent  developments were the introduction of twisted orbifold $K$-theory \cite{MR1993337} and the orbifold quantum cohomology \cite{MR2104605} on the topological side, and the investigation of gerbes \cite{MR2045679} and loop groupoids \cite{MR1950946} on the geometric side.

\subsubsection{}
Classically,  orbifolds are defined like manifolds as spaces which are locally homeomorphic to a quotient of an euclidean space by a finite group. Alternatively, orbifolds are represented by proper \'etale  smooth groupoids. Working with groupoid representations of orbifolds is like working with manifolds with a fixed atlas. In the modern coordinate invariant point of view an orbifold is a smooth stack in smooth manifolds which admits an orbifold atlas. By considering  orbifolds as objects in the $2$-catgeory of smooth stacks one makes the notion of morphisms\footnote{The right notion of a morphism between orbifolds is
 a representable morphism of stacks. This definition corresponds to the notion of a good morphism in the literature.} and other constructions like fibre products transparent.
The framework of stacks is most natural if one wants to include gerbes into the picture.

\subsubsection{}
If one replaces smooth manifolds by topological spaces then the corresponding analog of an orbifold is an orbispace. The goal of the present paper is to show that some geometric constructions on orbifolds are in fact topological concepts and extend to orbispaces. 

\subsubsection{}

The fixed point manifolds of the elements of the local automorphism groups of an orbifold $X$ can be assembled into a new orbifold $LX$ called the inertia or loop orbifold or the orbifold of twisted sectors. In the present paper we show that the loop orbifold can be characterized as the $2$-categorial (in the $2$-category of stacks) equalizer of the pair $(\id_X,\id_X)$ of the identity morphisms. The same definition applies to orbispaces in the topological context.  Since $2$-categorial equalizers always exist in $2$-categories of stacks it is clear that $LX$ exists as a stack. But it is not a priori clear that $LX$ is again an orbifold (or orbispace, resp.).  In the present paper we show that taking loop stacks preserves orbispaces. The corresponding result for orbifolds is well-known, but requires different, manifold specific arguments.

\subsubsection{}
A $U(1)$-banded gerbe $G\to X$ over an orbifold gives rise to a $U(1)$-principal bundle
$\tilde G\to LX$ over the loop orbifold of $X$. This bundle has a natural reduction of structure groups to the discrete $U(1)^\delta$. The traditional  way to construct this reduction is to choose a connection and curving on the gerbe $G\to X$. This geometric data induces a connection on $\tilde G\to LX$ which turns out to be flat by a calculation. The flat connection gives the reduction of structure groups, which turns out to be independent of the choices of geometric structures. The sheaf $\cL$ of locally constant sections of the associated flat line bundle $L\to LX$ is called an inner local system and plays an important role in the definition of the orbifold cohomology.

In  the present paper we give a topological construction of the reduction of the structure group of $\tilde G\to LX$ to $U(1)^\delta$ and of the sheaf $\cL$. Furthermore, we calculate its holonomy  in terms of the Dixmier-Douady class of the gerbe $G\to X$.

\subsubsection{}
The third concept which we generalize to the topological case is that of twisted delocalized orbifold cohomology. The usual definition in the smooth case is based on the de Rham complex of forms on $LX$ with coefficients in $L\to LX$. The differential of this complex involves the flat connection on $L$ and a closed three form on $X$ which represents the image of the  Dixmier-Douady class of the gerbe $f:G\to X$ in real cohomology. Let $f_L:G_L\to LX$ denote the pull-back of the gerbe via $LX\to X$. In the present paper we use the sheaf theory for smooth (or topological, resp.)  stacks
\cite{bss} in order to define the twisted delocalized orbifold cohomology as sheaf cohomology
$H^*(LX, \Tw_G(\cL))$, where $\Tw_G(\cL):=R(f_L)_*f_L^*(\cL)$. Our main result is, that in the smooth case the twisted delocalized cohomology according to this sheaf theoretic definition is isomorphic to the former construction using the de Rham complex. In addition to the fact that it works in the topological context our sheaf theoretic  definition of twisted delocalized orbifold cohomology has the advantage  that it is functorial in the gerbe $G\to X$.

\subsubsection{}

In the remaining parts of the introduction we give a detailed description of the results of the present paper and explain how they are related to the existing literature.

\subsection{A  description of the results}

\subsubsection{}

In the present paper we consider stacks in smooth manifolds or stacks in topological spaces.
 Our basic reference for stacks in these contexts is  \cite{heinloth}, but see also \cite{math.AG/0503247}, \cite{math.DG/0306176} and the recent \cite{math.DG/0605694}. A stack $X$ in smooth manifolds (topological spaces, resp.) is called a smooth stack (topological stack, resp.) if it admits an atlas $A\to X$.  The atlas is called an orbifold (orbispace, resp.) atlas if the smooth (topological, resp.) groupoid $A\times_XA\Rightarrow A$ is proper \'etale (very proper, \'etale and separated (see \ref{uiiuduiwed} for explanations)). An orbifold (orbispace, resp.) is a smooth (topological, resp.) stack which admits an orbifold (orbispace, resp.) atlas.

We refer to \cite{math.GT/0508550} for an introduction to orbispaces, and e.g. to \cite[Sec 2.]{MR2104605} for some basic information on orbifolds.

\subsubsection{}

In Subsection \ref{uiiuiqdww5} we review the notion of $2$-categorial limits.
The $2$-categorial equalizer of a pair of maps  is a special kind of limit.
We will see that equalizers exist in the $2$-catgeory of stacks on a site and in the two-category of groupoids in topological spaces.

The goal of Subsections \ref{ttzwtzetzwiu} and \ref{wsdwkjdwds} of the present paper is to place the construction of the loop orbifold $LX$ (or orbispace, resp.) into the framework of stacks in manifolds (topological spaces, resp.). 

We consider the orbifold (orbispace, resp.) $X$ as a stack and define its inertia stack $IX\to X$ as the $2$-categorial equalizer of the pair $(\id_X,\id_X)$.
The loop stack $LX$ is defined in an ad-hoc manner, see Definition \ref{uhdwc} and Remark \ref{equforstac}. We will see that it is canonically equivalent to $IX$. Though  Definition \ref{zwezuzu} of the $2$-categorial equalizer by a pull-back diagram is quite constructive we prefer to work with the simpler construction  $LX$ from now on. If $X$ is an orbifold (orbispace, resp.), then apriori $LX$ is a stack in smooth manifolds (topological spaces, resp.). 
In Lemma \ref{ll3}  (Lemma \ref{ll56}) we show that the loop stack of a topological stack (orbispace, resp.) is again a topological stack (orbispace). The main idea is to show that the existence of an (orbispace)  atlas of $X$ implies the existence of an (orbispace) atlas of $LX$. 

In the smooth case, the fact that the loop stack of an orbifold is again an orbifold is well known, see \cite[Lemma 3.1.1]{MR2104605} or \cite[Cor. 2.6.2]{MR1950946}.

\subsubsection{}

The loop orbifold  is also known as the orbifold of twisted sectors (compare \cite[Sec. 3.1]{MR2104605}) or inertia orbifold. It
plays an important role in the construction of the delocalized orbifold cohomology. The twisted sectors first appeared in connection with the orbifold index theorem
\cite{MR0474432}, \cite{MR641150}. In the framework of topological groupoids $G$ the corresponding object is called the inertia groupoid $\Lambda G$ which has been studied in detail in \cite{MR1950946}.  In order to keep our notation uniform in the present paper we will denote the inertia groupoid by $LG$ and call it loop groupoid\footnote{Note that the  loop groupoid $LG$ in
 \cite{MR1950946} is a much bigger object, and it is related with $\Lambda G$ by the equation $LG^\R\cong \Lambda G$  in the notation of \cite[Prop. 3.6.6]{MR1950946}.}.

\subsubsection{}\label{dkjsdsssss}

To a topological group $G$ we associate the classifying stack $\cB G:=[*/G]$
(see \cite[Example 1.5]{heinloth}). A $G$-principal bundle over a stack $X$ is by definition
a map $p:X\to \cB G$\footnote{Sometimes we will use a more sloppy language and say that $E\to X$ is a $G$-principal bundle, where $E\to X$ is defined by the pull-back
$$\xymatrix{E\ar[r]\ar[d]&{*}\ar[d]\\X\ar@{:>}[ur]\ar[r]^p&\cB G}\ .$$}.
Applying the loop functor and using the canonical isomorphism $L\cB G\cong [G/G]$ we get a map $Lp:LX\to [G/G]$. If $G$ is abelian, then this map lifts to a function
$h:LX\to G$. We are in particular interested in the case $G=U(1)$ and give various geometric and cohomological interpretations of this function.

In the present paper, ordinary cohomology of an orbispace $X$ is understood in the sense of \cite[Sec. 2.2]{math.GT/0508550}. Let $A\to X$ be an atlas and form the simplicial space $A^\cdot$ such that
$A^n:=A\times_X\dots\times_XA$ ($n+1$-factors). Here the fibre product is taken in stacks in topological spaces, but the stack $A^n$ is in fact  equivalent to a space since the map  $A\to X$ is representable. The cohomology of $X$ with integral coefficients is then  defined as
$$H^*(X;\Z):=H^*(|A^\cdot|;\Z)\ ,$$
where $|A^\cdot|$ denotes the realization of the simplicial space. Independence of the choice of the atlas has been shown in  \cite[Sec. 2.2]{math.GT/0508550} and \cite{MR2172499}\footnote{The  result in this paper is more general. The only condition on the atlas $A\to X$ is that the range and source maps of the groupoid $A\times_XA\Rightarrow A$ are topological submersions.}.   An alternative definition of the cohomology of $X$ could be based on the sheaf theory for orbifolds which will be discussed below. 
The group $H^2(X;\Z)$ classifies isomorphism classes of $U(1)$-principal bundles $p:E\to X$ (see \cite[Sec. 4.2]{math.GT/0508550} for this fact). 
\begin{enumerate}
\item
% For a finite group $\Gamma$ there is a natural identification of the group cohomology  $H^2(\Gamma;\Z)$ with the group of $U(1)$-characters $\hat \Gamma$ (see \ref{ejfwws}). If $\Gamma$ acts on a space $W$, then by $[W/\Gamma]$ we denote the quotient stack .
% In particular we can consider the action of $\Gamma$ on the one-point space $*$ and have the topological stack $[*/\Gamma]$.  Its loop stack  is
% $L[*/\Gamma]=[\Gamma/\Gamma]$, where $\Gamma$ acts on itself by conjugation. 
% 
% 
% Let $A:=*$ and consider the
% canonical map $A\to [*/\Gamma]$ which is an atlas of the stack $[*/\Gamma]$. The space
% $|A^\cdot|$  represents the homotopy type of the classifying space  $B\Gamma$ of $\Gamma$. 
% 
% Therefore, we  have
% $$H^*([*/\Gamma];\Z)\cong H^*(B\Gamma;\Z)\cong H^*(\Gamma;\Z)\ .$$
% In particular, we
% have an isomorphism
% $$\hat \Gamma \cong H^2(\Gamma;\Z)\cong H^2([*/\Gamma];\Z)\ .$$ Therefore a class $\chi\in H^2([*/\Gamma];\Z)$ gives rise to a function $\bar \chi:L[*/\Gamma]\to U(1)$,  where we view $U(1)$-valued characters of $\Gamma$ as functions  
% $L[*/\Gamma]=[\Gamma/\Gamma]\to U(1)$.
If $\Gamma$ is a finite group, then we have $H^2([*/\Gamma];\Z)\cong H^2(\Gamma;\Z)\cong H^1(\Gamma;U(1))\cong \hat \Gamma$, where $\hat \Gamma:=\Hom(\Gamma,U(1))$.  A class $\chi\in H^2([*/\Gamma];\Z)$ thus gives rise to  function $\bar \chi:L[*/\Gamma]\cong [\Gamma/\Gamma]\to U(1)$. This construction extends to general orbispaces $X$ and associates to each class
$\chi\in H^2(X;\Z)$ a function $\bar \chi:LX\to U(1)$.
A class $\chi\in H^2(X;\Z)$ also classifies a $U(1)$-principal bundle and therefore gives ries to a function
$h:LX\to U(1)$.
 We will show that
$\bar \chi=h$. This equality
has the following geometric interpretation. A point in the fibre of $LX\to X$ over $x\in X$ is an element $\gamma\in \Aut(x)$. The group $\Aut(x)$ acts on the fibre
$E_x:=p^{-1}(x)$. We write this as a left action. Then we  show in Lemma \ref{jujzsdh} that the functions $h, \bar \chi$ are both  characterized by 
$$\gamma e=e\bar\chi(\gamma)\ ,\quad \gamma e=eh(\gamma)$$
for all $e\in E_x$.  
\item $G$-principal bundles can be defined in terms of cocycles. We will give an interpretation of
the function $h$ in terms of the cocyle.
 \item A third cohomological interpretation uses the transgression  $\Tr:H^2(X;U(1))\to  H^1(LX;U(1))$ introduced in  \cite{math.AT/0605534}, \cite{math.AT/0307114}, \cite{math.KT/0604160}. 
 \end{enumerate}

\subsubsection{}

Let $f:G\to X$ be a topological gerbe with band $U(1)$ over an orbispace $X$.
The induced map $Lf:LG\to LX$ can be factored  canonically as $LG\stackrel{p}{\to} G_L\stackrel{f_L}{\to}LX$, where  $G_L:=LX\times_XG$. Here $f_L:G_L\to LX$ is a topological
gerbe with band $U(1)$, and $p:LG\to G_L$ is (the underlying map of) a $U(1)$-principal bundle.
The first main observation of Subsection \ref{twto} is that the bundle
$LG\to G_L$ descends canonically to a $U(1)$-principal bundle $\tilde G\to LX$.
The second result of this Subsection assumes that $X$ is an orbispace and asserts that $\tilde G\to LX$ has a canonical reduction of the structure group $\tilde G^\delta\to LX$ to  $U(1)^\delta$, the group $U(1)$ with the discrete topology.

The heuristic picture is as follows. 
Roughly speaking, a gerbe $G\to X$ over a topological stack $X$ associates to each point $x\in X$ a $U(1)$-central extension $1\to U(1)\to \widehat{\Aut(x)}\to\Aut(x)\to 1$ of the group of automorphisms $\Aut(x)$.
The fibre of the canonical map  $LX\to X$ over $x\in X$ is the automorphism group $\Aut(x)$. The $U(1)$-principal bundle $\tilde G\to LX$ restricts to $\widehat{\Aut(x)}\to \Aut(x)$ over $x\in X$. If $X$ is an orbispace, then finiteness of the groups $\Aut(x)$ provide a reduction of the structure group of this bundle to $U(1)^\delta$.  

 Let $L\to LX$ denote the  complex line bundle associated to $\tilde G^\delta\to LX$.  Since its structure group is discrete we can form the sheaf $\cL$ of locally constant sections of $L$.   

By \ref{awzegbdsajmnd} we have actually an extension
$$X\times U(1)^\delta\to \tilde G^\delta \to LX$$
of group stacks over $X$. The induced algebraic structures on $L\to LX$ turn this line bundle into an inner local system in the sense of \cite[Def. 2.1]{math.AG/0005299}, \cite[Def. 2.2.2]{MR2045679}.

\subsubsection{}
 
In the framework of  groupoids the construction of $\tilde G^\delta\to LX$ has been  previously given in   \cite[Thm. 6.4.2]{MR1950946} and \cite[Prop. 2.9]{math.KT/0505267} (with the exception of the reduction of the structure group to the discrete $U(1)^\delta$).
In the smooth case a reduction of the structure group of a line bundle from $U(1)$ to $U(1)^\delta$ is equivalent to a flat unitary connection. It has been observed in \cite[Lemma 5.0.1]{math.AT/0307114} and \cite[Prop. 3.9]{math.KT/0505267} that a connection on the gerbe $G\to X$ induces a flat connection on $L\to LX$.

Our original contribution here is to give a construction of this reduction of the structure group in purely topological terms. In addition to simplifications this extends the previous results to the topological case.

A twisted torsion in the language of
\cite{math.AG/0005299} is a class $\alpha\in H^2(\pi_1^{orbifold}(X),U(1))$, i.e. an isomorphism class of central $U(1)$-extensions
$$1\to U(1)\to \widehat{\pi_1^{orbifold}(X)}\to \pi_1^{orbifold}(X)\to 1\ .$$
The orbifold fundamental group $\pi_1^{orbifold}(X)$ is the automorphism group of the universal orbifold covering $Y\to X$. The map
$G_\alpha:=[Y/\widehat{\pi_1^{orbifold}(X)}]\to [Y/\pi_1^{orbifold}(X)]=X$
is a topological gerbe with band $U(1)$ over $X$.
 In 
\cite[Sec 4]{math.AG/0005299} or \cite[Example 2.2.2]{MR2045679} an inner local system $L_\alpha$ is associated directly to a twisted torsion $\alpha$. In the philosophy of the present paper we would consider $L_\alpha$ as the bundle associated to the gerbe $G_\alpha\to X$ via
the $U(1)^\delta$-bundle $\tilde G^\delta_\alpha\to LX$. 

The sheaf of locally constant sections $\cL$ of the line bundle $L$ (also called  inner local system) plays an importand role in the definition of twisted delocalized cohomology of an orbifold \cite{MR1993337}, \cite[Def. 2.2]{math.AG/0005299}\footnote{This is the cohomology of $LX$ with coefficients in $\cL$ with shifted grading. It is different from the gerbe-twisted delocalized cohomology}, \cite[Def. 3.10]{math.KT/0505267}.

\subsubsection{}

It is an interesting problem to calculate the holonomy of the bundle $\tilde G^\delta\to LX$ in terms of the Dixmier-Douady class $d\in H^3(X;\Z)$. We discuss this question in a typical case in Subsection \ref{cycl1}. Let $\pi:E\to X$ be a $U(1)$-principal bundle in orbispaces and $G\to E$ be a topological gerbe with band $U(1)$ and Dixmier-Douady class $d\in H^3(E;\Z)$. Let $\chi\in H^2(X;\Z)$ be the first Chern class of $E\to X$.
As explained in \ref{dkjsdsssss} we get  a function $\bar \chi:LX\to U(1)$.
Let $LX_1:=\bar\chi^{-1}(1)$. We will see that the canonical map $LE\to LX$ factorizes over $LX_1$, and that
$LE\to LX_1$ is again a $U(1)$-principal bundle.
The holonomy of the  bundle $\tilde G^\delta\to LE$ along the fibres of $LE\to LX_1$  can be considered as a function $$g:LX_1\to U(1)\ .$$ Our main result is the following description of this function. Note that $\pi:E\to X$ is an oriented fibre bundle.
We have an integration map $\pi_!:H^3(E;\Z)\to H^2(X;\Z)$. In particular we can form
$\pi_!(d)\in H^2(X;\Z)$ and the associated function $$\overline{\pi_!(d)}:LX\to U(1)\ .$$
In Proposition \ref{nhedvhjeada} we show the equality of functions
$$g=\overline{\pi_!(d)}_{|LX_1}\ .$$
In the smooth case (i.e. for orbifolds) holonomy questions could be addressed using Deligne cohomology. In fact, Deligne cohomology $H^*_{Del}(X)$ for orbifolds has been introduced in \cite{math.AT/0307114}. The choice of a connection on the gerbe $G$ leads to a lift of the Dixmier-Douady class $d\in H^3(X;\Z)$ of $G\to X$ to a Deligne cohomology class $d_{Del}\in H_{Del}^3(X)$ under the natural forgetful map $H^3_{Del}(X)\to H^3(X;\Z)$.  The transgression of $d_{Del}$  according to \cite[Thm. 6.0.2]{math.AT/0307114} is a class
$\Tr(d_{Del})\in H_{Del}^2(LX)$. Its integral $(L\pi)_!(\Tr(d_{Del}))\in H^1_{Del}(LX_1)$ should\footnote{We have not checked the details here. In this picture it is also not obvious that $(L\pi)_!(\Tr(d_{Del}))$ only depends on $d\in H^3(X;\Z)$, and not on the choice of its lift $d_{Del}\in H^3_{Del}(X)$.} give the function $g:LX_1\to U(1)$.

\subsubsection{}
\newcommand{\bLX}{\mathbf{LX}}
\newcommand{\bGL}{\mathbf{G_L}}
Section \ref{zuwlklklklkdee} of the present paper is devoted to twisted delocalized cohomology.
We are in particular interested in a version which is the target of the Chern character from twisted $K$-theory.  We refer to Subsection \ref{hhjjhjhjhdf} for a detailed introduction and a motivation of the particular definition of twisted delocalized cohomology. Our main original contribution in the present paper is a construction
of this cohomology in the framework of sheaf theory on topological stacks.
All previous definitions used the de Rham complex and are therefore tied to the orbifold case. 

To a topological stack (smooth stack, resp.) $X$ we associate a site $\bX$. 
The smooth case was discussed at length in \cite{bss}. So let us fix our conventions for the topological case here. A detailed account for the sheaf theory on topological stacks will be given in the paper \cite{bssf}. 

 An object of $\bX$ is a map $(\phi:U\to X)$ in stacks in topological spaces, where $U$ is a topological space (or more precisely stack which is equivalent to a space), and $\phi$ is a representable map  which admits local sections\footnote{Note that $\bX$ must be small. A precise definition would either involve universes or a cardinality restriction.}. 
The morphisms in $\bX$ are commutative diagrams
$$\xymatrix{U\ar[rr]\ar[dr]&\ar@{:>}[d]&V\ar[dl]\\&X&}$$
consisting of a morphism $U\to V$ and a $2$-morphism. A family $(U_i\to U)_{i\in I}$ of morphisms in $\bX$ is a covering family if all maps $U_i\to U$ admit local sections and the induced map $\sqcup_{i\in I} U_i\to U$ is surjective. To the site $\bX$ we can associate the category of sheaves $\Sh\bX$ of sets and the abelian category $\Sh_\Ab\bX$ of sheaves of abelian groups.

A map between topological (resp. smooth) stacks $f:X\to Y$ induces an adjoint pair of functors $$f^*:\Sh\bY\Leftrightarrow \Sh\bX: f_*$$ relating the categories of sheaves on these sites. In the smooth case the construction of this adjoint pair was given by \cite[Sec. 2.1]{bss}. The construction in the case of topological stacks is very similar, see \cite{bssf}.

The restriction $f_*:\Sh_\Ab\bX\to \Sh_\Ab\bY$ of $f_*$ to abelian sheaves is left-exact and admits a right-derived functor $$Rf_*:D^+(\Sh_\Ab\bX)\to D^+(\Sh_\Ab\bY)$$ between the lower-bounded derived categories.

 Let $G\to X$ be a topological (smooth, resp.) gerbe with band $U(1)$ on an orbispace (resp. orbifold)  $X$. It gives rise to the $U(1)^\delta$-principal bundle $\tilde G^\delta\to LX$ and an associated locally constant sheaf $\cL$ of $\C$-vector spaces on the site $\bLX$. 
In Subsection \ref{fdifdfdmhjjhhjjh} we define the $G$-twisted delocalized cohomology
of $X$ by
\begin{equation}\label{iwdhwkjwkwk}
H^*_{deloc}(X;G):=H^*(\ev\circ Rp_*\circ f_L^*(\cL))\ .
\end{equation}
The notation is explained by means of  the following diagram
$$\xymatrix{{*}&G_L\ar[l]^p\ar[d]^{f_L}\ar[r]&G\ar[d]^f\\&LX\ar@{:>}[ur]\ar[r]&X}\ ,$$
where the square is $2$-cartesian, i.e.  $f_L:G_L\to LX$ is the pull-back of the gerbe $f:G\to X$ via the canonical map $LX\to X$,  and the map $p:G_L\to *$ is the canonical projection to the point.  
Since $\Site(*)$ is the big site of the point, i.e. the category of all non-empty topological spaces we need the evaluation $\ev:D^+(\Sh_\Ab\Site(*))\to D^+(\Ab)$ at the object $(*\to *)\in \Site(*)$. The functoriality of this cohomology in the data $G\to X$ is studied in Lemma \ref{zuzuzwdwdwdw}.

Our main result is the comparison of this sheaf-theoretic definition of $G$-twisted delocalized cohomology with the previous de Rham model
\cite[Def. 3.10]{math.KT/0505267} in the case of orbifolds.

\subsubsection{}

We now explain the de Rham model for the twisted delocalized cohomology.
Let $X$ be an orbifold and $G\to X$ be a smooth gerbe with band $U(1)$.
In this case we can define three versions of twisted delocalized de Rham cohomology. The $2$-periodic twisted delocalized cohomology is the correct target of the Chern character and will be defined in \ref{qwekjdqjdqqdd}.
The sheaf theoretic cohomology  (\ref{iwdhwkjwkwk}) is not $2$-periodic. In the following we describe its appropriate de Rham model.
%The one can define the $G$-twisted delocalized $2$-periodic cohomology %of $X$  using the de Rham complex \cite[Def. 3.10]{math.KT/0505267},see %also \ref{qwekjdqjdqqdd}.
We choose a closed three-form $\lambda\in \Omega^3(LX)$ which represents the image of the Dixmir-Douady class of $G_L\to LX$ in real cohomology. Then we define a sheaf $\Omega_{LX}[[z]]_\lambda\in C^+(\Sh_\Ab \bLX)$ of complexes which associates to each object $(\phi:U\to LX)\in \bLX$ the complex $(\Omega(U)[[z]],d_\lambda)$,
where $(\Omega(U),d_{dR})$ is the de Rham complex of the smooth manifold $U$, $z$ is a formal variable of degree $2$, and $d_\lambda=d_{dR}+\frac{d}{dz} \phi^*\lambda$. 
Let $\Omega(LX;\cL)[[z]]_\lambda:=\Gamma_{LX}(\Omega_{LX}[[z]]_\lambda\otimes \cL)$ denote the complex of global sections (see \ref{uieuhwuefuw} for the definition of global sections) of the tensor product of sheaves $\Omega_{LX}[[z]]_\lambda\otimes \cL$.
Its cohomology is the twisted delocalized de Rham cohomology \begin{equation}\label{uzwdwjkkj}
H^*_{dR,deloc}(X,(G,\lambda)):=H^*(\Omega(LX;\cL)[[z]]_\lambda)
\end{equation}
(see \ref{def36}).

\subsubsection{}

The twisted delocalized de Rham cohomology defined in
\cite[Def. 3.10]{math.KT/0505267} is related  to the definition of the present paper by a duality. Let  us first recall
the definition \cite[Def. 3.10]{math.KT/0505267}. Let $u$ be a formal variable of degree $-2$ and define the complex of sheaves  $\Omega_{LX}((u))$ which associates to $(\phi:U\to LX)$ the space of formal Laurent series of forms
$\Omega(U)((u))_\lambda$ with the differential $d^\prime_\lambda:=d_{dR}-u  \imath\phi^*\lambda$.
The twisted cohomology in \cite[Def. 3.10]{math.KT/0505267} is
 the cohomology of the complex of compactly supported global sections $\Omega(LX;\cL)_{comp}((u))_\lambda$\footnote{Here we use the freedom of rescaling $\lambda$ by non-zero factors as explained in \cite[Rem. 3.11(1)]{math.KT/0505267}.} of
$\Omega_{LX}((u))_\lambda\otimes \cL$. 
Note that the multiplication by $u$ induces an isomorphism of complexes
which makes the cohomology of \cite[Def. 3.10]{math.KT/0505267} two-periodic.

We define the pairing (using the hermitean structure of $\cL$) 
\begin{equation}\label{hdjqdwq45}
<\dots,\dots>:\Omega(LX;\cL)_{comp}((u))_\lambda\otimes \Omega(LX;\cL)[[z]]_\lambda\to \C
\end{equation}
by $$<u^n\omega,z^m\alpha^l>=\delta_{m,n}m!\int_{LX}\omega\wedge \alpha\ ,$$
where $\omega\in  \Omega(LX;\cL)_{comp}$ and $\alpha\in \Omega(LX;\cL)$.
One easily checks that
$$<d_\lambda^\prime \omega,\alpha>=(-1)^{|\omega|+1}<\omega,d_\lambda \alpha>\ .$$
The pairing (\ref{hdjqdwq45}) induces an embedding
of $\Omega(LX;\cL)[[z]]_\lambda$ into the dual complex of $\Omega(LX;\cL)_{comp}((u))_\lambda$.

\subsubsection{}

Let us now explain the relation between (\ref{uzwdwjkkj})  and the $2$-periodic version  \ref{qwekjdqjdqqdd}. Note that
the complex of sheaves $\Omega_{LX}[[z]]_\lambda$
admits an action of the operation $T:=\frac{d}{dz}$ of degree $-2$.
We consider the system
$$\cS: \Omega_{LX}[[z]]_\lambda\stackrel{T}{\leftarrow} \Omega_{LX}[[z]]_\lambda[2]\stackrel{T}{\leftarrow} \Omega_{LX}[[z]]_\lambda[4]\stackrel{T}{\leftarrow}\dots
$$
in the category $C(\Sh_\Ab\bLX)$ of unbounded complexes.
The discussion of \cite[1.3.23]{bss} can be subsumed in the assertion that
$\Gamma_{LX}(\lim \cS\otimes \cL)$ is exactly the periodic complex (\ref{perzwzgehjjhb}).

\subsubsection{}

Our basic result, Theorem \ref{main1},  is an extension of \cite[Thm. 1.1]{bss} from smooth manifolds to orbifolds.
It asserts that there is an isomorphism
\begin{equation}\label{wejkjkkjw}R(f_L)_*( \R_{\bGL})\cong \Omega_{LX}[[z]]_\lambda\end{equation}
 in the dervied category $D^+(\Sh_\Ab\bLX)$. 
This isomorphism is not canonical and depends on the choice of a connection on the gerbe $G\to X$.
As a consequence of  (\ref{wejkjkkjw}) we get in Theorem \ref{main2} the non-canonical isomorphism 
$$H^*_{dR,deloc}(X;(G,\lambda))\cong H^*_{deloc}(X;G)\ .$$

\subsubsection{}
\newcommand{\holim}{{\tt holim}}
The main goal of the forthcoming paper \cite{bssm} will be  a sheaf theoretic  construction of $2$-periodic twisted delocalized cohomology. The idea is to define an analog $T$ of the operation $\frac{d}{dz}$ on the left-hand side of the derived category isomorphism  (\ref{wejkjkkjw}). In analogy with the de Rham model  we then will consider the system
$$\cT: R(f_L)_*(\R_{\bGL})\stackrel{T}{\leftarrow} R(f_L)_*(\R_{\bGL})[2]\stackrel{T}{\leftarrow}R(f_L)_*( \R_{\bGL})[4]\stackrel{T}{\leftarrow}\dots
$$
in $D(\Sh_\Ab \bLX)$.
The sheaf-theoretic version
of periodic delocalized twisted cohomology will be  defined as
$$H^*(\ev\circ Rp_*( \holim \cT \otimes \cL))\ .$$
In order to make this rough idea precise we must solve various problems, in particular
\begin{enumerate}
\item The homotopy limit $\holim\cT$ of the digragram $\cT$ in the derived category is only well-defined up to non-canonical isomorphism. In order to define a functorial periodic cohomology we must work hard to construct a much more concrete version of the system $\cT$.
\item The push-forward $Rp_*(\colim \cT\otimes \cL)$ is not a standard derived functor since it acts between unbounded derived categories. We use a model category approach in order to construct functors like $Rp_*$.
\end{enumerate}
The main application and technical tool in \cite{bssm}  will be   $T$-duality. 
The results of Subsections \ref{hdjjhsjcoo} and \ref{cycl1} of the present paper will be  needed in \cite{bssm} in a crucial way.

\subsection{Motivation of the definition of twisted delocalized cohomology}\label{hhjjhjhjhdf}

\subsubsection{}

In the present subsection we motivate the definition of twisted delocalized cohomology as the correct target for the Chern character from twisted $K$-theory.

It is a well-known fact that the Chern character $\ch:K(X)\rightarrow H(X;\Q)$ from the complex $K$-theory of a space $X$ to the rational cohomology of $X$ induces an isomorphism $K(X)\otimes_\Z \Q\stackrel{\sim}{\rightarrow} H(X;\Q)$
(we consider both sides as $\Z/2\Z$-graded groups)

\subsubsection{} 

Complex $K$-theory and rational cohomology both have equivariant generalizations. Every generalized cohomology $E$ theory has the Borel extension. If $X$ is a $G$-space, then the Borel extension of $E$ to $G$-spaces associates to $X$ the group
$E^{Borel}_G(X):=E(EG\times_G X)$. Here $EG$ is a universal  space for $G$, i.e. a contractible space on which $G$ acts freely. The Chern character induces an equivariant Chern character $\ch_G:K_G^{Borel}(X)\rightarrow H_G^{Borel}(X;\Q)$ which gives again a rational isomorphism.

\subsubsection{}

The interesting equivariant extension of $K$-theory is not the Borel extension but the extension due to Atiyah-Segal based on equivariant vector bundles \cite{MR0259946}. It will be denoted by $K_G(X)$. In order to see the difference between $K_G^{Borel}$ and $K_G$ consider the simple example of finite group $G$ acting trivially on the point $*$. The equivariant Atiyah-Segal $K$-theory is isomorphic to the representation ring $R(G)$ of $G$. In \cite{MR0148722} is was shown that $K^{Borel}_G(*)$ is isomorphic to the completion $\widehat{R(G)}_{I}$ of the representation ring at the dimension ideal $I$, which is defined as the kernel of the homomorphism $\dim:R(G)\to \Z$.

\subsubsection{}\label{dhswhsd}

It is not true that the Atiyah-Segal equivariant $K$-theory is rationally isomorphic to the Borel extension of rational cohomology. In the case of discrete groups and proper  actions the appropriate target of the Chern character was found in \cite{MR928402}. It will be called the delocalized cohomology in this paper. Let $G$ be a discrete group which acts properly on a space $X$. Then we define a new proper $G$-space (sometimes called the Brylinski space)
$$\Lambda X:=\bigsqcup_{g\in G} X^g\ ,$$ where $X^g\subset X$ is the subspace of fixed points of $g$. The action of $h\in G$ on $\Lambda X$ maps $x\in X^g$ to $hx\in X^{hgh^{-1}}$. The delocalized cohomology of the $G$-space $X$ is the cohomology of the quotient $\Lambda X/G$.

\subsubsection{}

A $G$-space $X$ gives rise to a topological quotient stack $[X/G]$. 
% We refer to \cite{heinloth} and \cite{math.AG/0503247}
% for the language of topological stacks. 
If $G$ is a discrete group which acts properly on $X$, then the quotient $[X/G]$ is an example of an orbispace (the topological variant of an orbifold). But  not every orbispace can be represented in this form. We refer to \cite{math.GT/0508550} for the description of the category of orbispaces.
The stack $[\Lambda X/G]$ has a description in the language of topological stacks. 
If $Z$ is a topological stack, then we  define its loop stack $LZ$ (see \ref{uhdwc} and \ref{equforstac})\footnote{In the present paper we use the name loop stack. In the literature it is also known under the name inertia stack} such that
$$L[X/G]=[\Lambda X/G]$$ for a discrete group acting properly on a space $X$.

\subsubsection{}

If $G$ is a discrete group which acts properly on a space $X$, then the quotient $X/G$ is a reasonable topological space. It is the coarse moduli space of the orbispace $[X/G]$.
The definition of the coarse moduli space extends to arbitrary orbispaces.
The coarse moduli space of the orbispace $Z$ will be denoted by $|Z |$. If 
$Z^1\Rightarrow Z^0$ is a presentation of the orbispace by a proper \'etale groupoid, then $|Z|=Z^0/Z^1$.

The rational  cohomology of an orbispace $Z$ is the cohomology of its coarse moduli space $|Z|$. Therefore we can define the delocalized cohomology of an orbispace as the cohomology of $|LZ|$. This generalizes the definition of the delocalized cohomology from global quotient orbispaces to general orbispaces.

Note that this is not quite the definition of delocalized cohomology which we are going to use in the main part of the paper but sufficient for the present 
discussion. Later we prefer a sheaf-theoretic definition
of the delocalized cohomology. 

\subsubsection{}

Delocalized cohomology for orbifolds appeared in connection with the index theorem for orbifolds \cite{MR641150}.
In a completely different context of quantum cohomology for orbifolds it was constructed in \cite{MR2104605}, \cite{MR1941583}. Note that the grading used in \cite{MR2104605} is different from the grading in the present paper.

\subsubsection{}\label{lambdat}

A different generalization of $K$-theory is twisted $K$-theory (see \cite{MR2172633}). The search for the target of an appropriate Chern character lead to the definition of $2$-periodic twisted de Rham cohomology\footnote{This could also be reversed. The equations for fields associated to $D$-branes in string theory with $B$-field backgroup can be expressed in terms of the twisted de Rham differential. In this history twisted $K$-theory was found as a cohomology theory with a (Chern character) map to twisted de Rham cohomology giving the integrality lattice  of $D$-brane charges \cite{Minasian:1997mm}, \cite{Witten:1998cd}.}. Usually it is defined on smooth manifolds $X$. Given a closed three-form $\lambda\in \Omega^3(X)$ twisted de Rham cohomology is the cohomology of the complex  
\begin{equation}\label{zezdwe7we7d87we989}
\dots\stackrel{d_\lambda}{\rightarrow} \Omega^{even}(X) \stackrel{d_\lambda}{\rightarrow} \Omega^{odd}(X)\stackrel{d_\lambda}{\rightarrow} \Omega^{even}(X)\stackrel{d_\lambda}{\rightarrow}\dots\ ,
\end{equation}
where $d_\lambda:=d_{dR}+\lambda$. 

\subsubsection{}

A Chern character for twisted $K$-theory with values in $\lambda$-twisted de Rham cohomology was constructed in  \cite{MR1911247}, \cite{MR1977885},  and \cite{math.KT/0510674}. The twist of $K$-theory is classified by a class $\lambda_\Z\in H^3(X;\Z)$. The closed form $\lambda\in \Omega^3(X)$ should represent the image of $\lambda_\Z$ in real cohomology.
It was shown that this $\Z/2\Z$-graded cohomology theory is again isomorphic to  twisted $K$-theory tensored with $\R$.

\subsubsection{}

Twisted $K$-theory on orbifolds has first been considered in \cite{MR1993337}. In this paper the twist was given by a so-called inner local system of twisted torsion. The natural object to be used to twist complex $K$-teory should a gerbe  $G\to X$ with band $S^1$. Gerbe twisted $K$-theory for orbifolds was discussed in \cite{MR2045679}. For general local quotient stacks it was defined in \cite{math.AT/0312155}.  Using topological groupoids in order to represent stacks a very general definition of twisted $K$-theory was given in \cite{MR2119241}.

\subsubsection{}\label{naiiv}

The result of \cite{MR928402} in the case of global quotient orbispaces obtained from proper actions of discrete groups shows that the correct target of the Chern character has to take the topology of the fixed point sets into account.
Thus the target of the Chern character from twisted $K$-theory of an orbifold should be a delocalized version of twisted de Rham cohomology.
If $X$ is an orbifold, then $LX$ is again an orbifold. In particular we can consider differential forms on $LX$. Given a three-form $\lambda\in \Omega^3(LX)$ we can define the twisted delocalized de Rham cohomology as the cohomology of the complex
\begin{equation}\label{pefffrzwzgehjjhb}
\dots\stackrel{d_\lambda}{\rightarrow} \Omega^{even}(LX) \stackrel{d_\lambda}{\rightarrow} \Omega^{odd}(LX)\stackrel{d_\lambda}{\rightarrow} \Omega^{even}(LX)\stackrel{d_\lambda}{\rightarrow}\dots\ .
\end{equation}
It turned out that this cohomology is not the correct target of the Chern character. 
This has already been observed at the end of \cite{MR1993337}. 

\subsubsection{}\label{qwekjdqjdqqdd}

Let $(L,\nabla^L)$ be the flat  complex line bundle  associated to $\tilde G^\delta\to LX$. We let $\Omega(LX;L)$ denote the differential forms with values in $L$, and $d^L$ be the differential induced by $d_{dR}$ and the flat connection
$\nabla^L$. We let $\lambda\in \Omega^3(LX)$ be a closed three from which represents the image of the Dixmir-Douady class $\lambda_\Z\in H^3(LX;\Z)$ of the gerbe $G_L\to LX$ in real cohomology. We set $d_\lambda^L:=d^L+\lambda$.
The correct target of the Chern character on $G$-twisted $K$ of the orbifold $X$ is the $2$-periodic cohomology of the complex
\begin{equation}\label{perzwzgehjjhb}\dots\to \Omega^{ev}(LX;L)\stackrel{d^L_\lambda}{\to}\Omega^{odd}(LX;L)\stackrel{d^L_\lambda}{\to}\Omega^{ev}(LX;L)\to \dots\ .
\end{equation}
This Chern character was constructed in 
 \cite{math.KT/0505267}.

\section{Inertia}\label{seczu2s}

\newcommand{\caC}{{\cC}}
\newcommand{\caD}{{\cD}}
\newcommand{\uHom}{{\underline{\tt Hom}}}
\newcommand{\uhom}{{\underline{\tt hom}}}
\newcommand{\Ob}{{\tt Ob}}
\newcommand{\Cat}{{\tt Cat}}
\newcommand{\PSt}{{\tt PSt}}
\newcommand{\St}{{\tt St}}
\newcommand{\groupoids}{{\tt gpd}}
\newcommand{\gpd}{\groupoids}
\newcommand{\caS}{{\mathcal{S}}}

\subsection{$2$-limits in $2$-categories}\label{uiiuiqdww5}

\subsubsection{}

In the present paper we consider stacks on some site or groupoids in some ambient category like topological spaces or manifolds. A common feature of these constructs is that they are objects in a $2$-category. Of particular importance for the present paper is the notion of a $2$-limit. The goal of this Subsection is to explain this notion.

\subsubsection{}

By a {\em $2$-category} we always mean a strict $2$-category. In our main examples of $2$-categories have the property that all
$2$-morphisms are isomorphisms, but in the present subsection do not assume this.
For objects $a$ and $b$ of a $2$-category
we denote by $\Hom_\cC(a,b)$ the Hom-category from $a$ to $b$ (we will often omit the subscript and write $\Hom(a,b)$).
We will write the objects as straight arrows $a\to b$, and the morphisms between 
two arrows $f,g:a\to b$ as $f\leadsto g$.

\subsubsection{} \label{2-functors}
By a $2$-functor we always mean a
pseudo-$2$-functor, as explained for example in
\cite[Definition 1.4.2]{MR1650134}. By a {\em strict} $2$-functor we mean such a functor where
all unit and composition $2$-isomorphisms are identities.

\subsubsection{}\label{uiduwed}

Let $\caC$ be a $2$-category. For any $X \in \Ob \caC$ we denote by $\caC/X$ the
over $2$-category 
\begin{itemize}
\item with objects the $1$-arrows $A \to X$,
\item  whose $1$-morphisms are 
triangles filled in with a $2$-morphism
$$
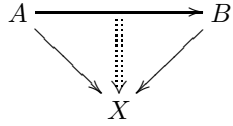
\xymatrix{A\ar[rr]\ar[dr]&\ar@{:>}[d]&B\ar[dl]\\&X&}$$, \item and where $2$-morphisms are the ones of $\caC$ making the
natural diagram commutative.
\end{itemize} 
There is a version of this construction
for a $2$-functor $\caD \to \caC$ and an object $X$ of $\caC$, denoted
$\caD/X$. Note that if $\caD$ is a $1$-category then so is $\caD/X$.

\subsubsection{}

Let $\caC$ be a $2$-category and $D$ a small category.
Let  $F,G : D \to \caC$ be two
$2$-functors. A {\em natural $2$-transformation} $\varphi$ from $F$ to $G$
is an assignment of a $1$-morphism $\varphi(a): F(a) \to G(a)$ for any object $a$ of $D$
and a $2$-isomorphism $\varphi (f)$ for any $f: a \to b$ in $D$ filling
in the square
$$\xymatrix{F(a) \ar[r]^{\varphi( a)} \ar[d]^{F (f)} & G(a) \ar[d]^{G (f)} \\
F(b)\ar@{:>}[ur]^{\phi(f)} \ar[r]^{\varphi( b)} & G(b)} \text{,}$$
satisfying the obvious  compatibility for compositions of maps in $D$.

Let $\varphi, \psi: F \to G$ be two natural $2$-transformations.
A {\em modification} $t$ from $\varphi$ to $\psi$
consists of a $2$-morphism $t(a): \varphi (a) \leadsto \psi (a)$ for any object
$a$ of $D$ satisfying an again obvious  compatibility with the $\varphi(f)$ and
$\psi (f)$ for any map $f$ in $D$.

With these definitions the $2$-functors, the natural $2$-transformations
and the modifications form a $2$-category.

For $F,G$ as above we denote by
$\Hom_{\cC^\cD}(F,G)$ the corresponding category of natural transformations from $F$ to $G$.

\subsubsection{}

For an object $c$ of $\caC$ we denote by $\underline{D}_c$ the
constant diagram on $c$, i.e. the (strict) $2$-functor from $D$ to $\caC$
sending all objects to $c$ and all morphisms to the identity on $c$.

\begin{ddd}
Let $F: D \to \caC$ be a $2$-functor. A {\em $2$-limit} of $F$
is an object $c$ of $\caC$ together with a natural $2$-transformation
$\varphi: \underline{D}_c \to F$ such that for any object $T$ of $D$
the functor
$$\Hom_{\cC}(T,c) \to \Hom_{\cC^\cD}(\underline{D}_T,F)$$
given by composition with $\varphi$ is an equivalence
of categories.
\end{ddd}

The constant diagram functor $c\mapsto \underline{D}_c$ is a $2$-functor $\cC\to \cC^D$.
Note that  $F\in \cC^D$. Using $\underline{D}$ we form  the over $2$-category $\cC/F$.
By definition a $2$-limit $(c,\phi)$ of $F$ is an object of $\cC/F$.

For example a $2$-final object of $\caC$ is an object $c$ such that
for all objects $T$ of $\caC$ the projection from $\Hom(T,c)$ to the
point category is an equivalence.

\begin{lem} \label{2-final-in-overcat}
Let $u: \caC \to \caD$ be a $2$-functor between $2$-categories,
$X$ an object of $\caD$. Let $c$, $f: u(c) \to X$ be an object of $\caC/X$.
Then if the functor
$$\Hom_\cC(T,c) \to \Hom_\cD(u(T),X)$$
is an equivalence for all objects $T$ of $\caC$ the object
$(c,f)$ is $2$-final in $\caC/X$. If the $2$-morphisms in $\caD$ and $\caC$ are all
$2$-isomorphisms the converse holds.
\end{lem}
\begin{proof}
Let $(c',f') \in \caC/X$ be another object. Then there is a canonical
$2$-cartesian square
$$\xymatrix{\Hom_{\caC/X}((c',f'),(c,f)) \ar[r] \ar[d] & \Hom_{\caC}(c',c) \ar[d] \\
\mathrm{pt} \ar[r]^{f'} & \Hom_{\caD}(u(c'),X)}$$
in $\Cat$. Hence the first statement follows.
The second statement follows from the fact that a map $\varphi: A \to B$ between groupoids
is an equivalence if and only if all ($2$-categorical) fibers over objects of $B$ are contractible.
\hB\end{proof}

\subsubsection{}

An equivalence between two objects $c$ and $d$ of $\caC$
are $1$-arrows $f: c \to d$ and $g: d \to c$ together with $2$-isomorphisms
$\varphi: \id_c \leadsto g \circ f$ and $\psi: \id_d \leadsto f \circ g$ satisfying
the triangular identities as for units and counits of adjunctions.

\subsubsection{} \label{2-final}
As particular case consider two $2$-final objects $c,c'$ in
a $2$-category $\caD$. Then there is an equivalence between
$c$ and $c'$ which is unique up to unique $2$-isomorphism.

\subsubsection{}

\begin{lem}  \label{limits-equivalent}
If an object $(c, \varphi) \in \caC/F$ is a $2$-limit of $F$ then it
is $2$-final in $\caC/F$. If all $2$-morphisms in $\caC$ are $2$-isomorphisms
or if $\caC$ has all small $2$-limits
then the converse is true.
Any two choices of $2$-limits are equivalent in $\caC/F$, unique up to unique $2$-isomorphism,
in particular the underlying objects in $\caC$ are (canonically) equivalent.
\end{lem}
\begin{proof}
The first statement follows from Lemma \ref{2-final-in-overcat}.
The second statement under the assumption on the $2$-morphisms also
follows from that Lemma, under the completeness assumption it follows
from the first statement and the uniqueness (up to unique isomorphism) of $2$-final objects \ref{2-final}.
The third statement is also \ref{2-final}.
\hB\end{proof}

\subsubsection{}

\begin{lem}
In $\Cat$, the $2$-category of small categories, small $2$-limits exist.
\end{lem}
\proof
The usual construction gives a {\em preferred model}: For a $2$-functor $F: D \to \Cat$
define $c$ to be the category whose objects are collections of objects $x_a \in F(a)$ for
any object $a$ of $D$ together with isomorphisms $\varphi_f: (Ff)(x_a) \to x_b$
for any map $f: a \to b$ in $D$ satisfying a compatibility condition for
compositions of maps in $D$, and whose morphisms from $(x_a)$
to $(y_a)$ are compatible
systems of morphisms $x_a \to y_a$.
The transformation $\underline{D}_c \to F$ induced by projections exhibits
$c$ as a $2$-limit of $F$.\hB

\subsubsection{}

Let us consider for example the category 
\begin{equation}\label{kkwkehdwe}
D:=\quad\xymatrix{&b\ar[d]\\a\ar[r]&c}\ .
\end{equation}
A functor $F:D\to \cC$ is a diagram
\begin{equation}\label{iiuwedwnn}
\xymatrix{&B\ar[d]^v\\A\ar[r]^u&C}\ .
\end{equation}
Usually a $2$-categorical fiber product of $F$
is a diagram 
\begin{equation}\label{wuiefiudwe}
\xymatrix{A\times_CB\ar[d]\ar[r]&B\ar[d]^v\\A\ar@{:>}[ur]^\psi\ar[r]^u&C}\
\end{equation}
fulfilling some natural properties. Such a diagram gives
in two natural ways an object in $\cC/F$ (by requiring the
map $A\times_CB \to C$ be one of the two possible compositions),
and it is easily checked that the usual properties are
equivalent to this object being a $2$-limit.

If these properties are fulfilled we call a diagram as above {\em $2$-cartesian}.

\subsubsection{}\label{zuzuuzewd}

Assume that $\cC=\Cat$. A model of $A\times_CB$ is then the category whose objects are triples $(a,b,\gamma)$, where $a\in \Ob(A)$, $b\in \Ob(b)$ and $\gamma:u(a)\to v(b)$. A morphism $(a,b,\gamma)\to (a^\prime,b^\prime,\gamma^\prime)$ is a pair $(f:a\to a^\prime,g:b\to b^\prime)$ such that $\gamma^\prime\circ u(f)=v(f)\circ \gamma$. 
The $2$-morphism in (\ref{wuiefiudwe}) is given by
$\psi(a,b,\gamma):=\gamma$.
 
We see in particular that $2$-categorical fiber products in $\Cat$ are $2$-limits.

We call  any diagram as \ref{wuiefiudwe} a {\em standard model}
of the fiber product in a $2$-category $\cC$ if for any
object $T$ the functor $\Hom(T,\_)$ produces a diagram
which is {\em isomorphic} (with respect to an obvious map)
to the model in $\Cat$ from above.
Note that this is not the preferred model.

\subsubsection{}

Like ordinary limits $2$-categorial limits are characterized by a universal property for Hom-categories.
\begin{lem} \label{hom-from-const-as-lim}
Let $F: D \to \caC$ be a $2$-functor, $(c,\varphi) \in \cC/F$
a $2$-limit of $F$ and $T$ and object of $\caC$.
Consider the $2$-functor $H:D\to \Cat$ given by
$a \mapsto \Hom_{\cC}(T,F(a))$.
Then the natural map $\underline{D}_ {\Hom_\cC(T,c)}\to H$ is a $2$-limit of the functor $H$.
\end{lem}
\begin{proof}
In fact $\Hom_{\cC^\cD}(\underline{D}_c,F)$ is naturally isomorphic to
the preferred model of the $2$-limit of the diagram
$a \mapsto \Hom_{\cC}(T,F(a))$.
\hB\end{proof}
 Lemma \ref{hom-from-const-as-lim} implies an equivalence of categories
$$\Hom_\caC(T,2-\lim_{a\in D}F(a))\cong  2-\lim_{a\in D} \Hom_{\caC}(T,F(a))\ ,$$
where the left $2$-limit is taken in $\caC$, and the right $2$-limit is taken in $\Cat$.

\subsubsection{}

Let $C$ be another small category and suppose given a $2$-functor
$F: C \times D \to \caC$. For simplicity suppose that $\caC$ has all
small $2$-limits.

\begin{prop}\label{hjhjhjhqa}
Let the notation be as above. The assignment
$$a \mapsto 2-\lim_{b \in D} F(a,b)$$ can be made into a $2$-functor
$K: C \to \caC$, and two such choices are canonically equivalent.
Moreover the $2$-limit of $K$ is canonically equivalent to
the $2$-limit of $F$.
\end{prop}
\proof
The first assertion is a consequence of Lemma \ref{limits-equivalent}.
We sketch  the proof of the second statement.
By Lemma \ref{hom-from-const-as-lim} we are reduced to prove
the statement in $\Cat$. But taking everywhere preferred models
produces isomorphic models of the two $2$-limits in question. 
\hB

\subsubsection{}\label{ifext}

We will assume that  $\cC$ has a final object and admits standard models (see \ref{zuzuuzewd}) of all
 $2$-categorical fiber products.
The absolute product $\times$ is understood as a standard model of the fiber product over the final object.
Consider a pair of maps 
$$
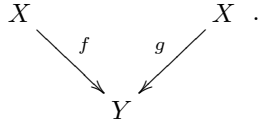
\xymatrix{X\ar[dr]^f&&X\ar[dl]_g\\&Y&}\ .$$

\begin{ddd}\label{zwezuzu}
The equalizer $E(f,g)$ of the pair of maps $f,g:X\to Y$ is defined as a standard model of the  $2$-categorical fiber product
$$
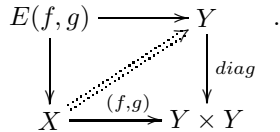
\xymatrix{E(f,g)\ar[d]\ar[r]&Y\ar[d]^{diag}\\
X\ar@{:>}[ur]\ar[r]^{(f,g)}&Y\times Y}\ .$$
\end{ddd}

Note that on Hom-categories this definition yields in fact the
preferred model of the equalizer diagram.

\begin{ddd}\label{iooiwed}
We define the inertia object of $X$ as the equalizer $IX:=E(\id_X,\id_X)$.
\end{ddd}

\subsubsection{}

We say that a $2$-category is $2$-complete if it admits a small $2$-limits. There is an analogous notion of a $2$-colimit, and the category is called $2$-cocomplete if all small $2$-colimits exist. The category is called $2$-bicomplete if it is $2$-complete and $2$-cocomplete.

The $2$-category
of small groupoids $\gpd$ is $2$-bicomplete as well as bicomplete as
a category. The same holds for the $2$-category $\PSt I$
of prestacks on a small category $I$, which is by definition the
$2$-category of $2$-functors  $\gpd^{I^{op}}$.
The $2$-category of stacks $\St \bS$ on a
small site $\bS$ is $2$-bicomplete.

\subsubsection{}
 
We consider the $2$-category $\groupoids(\cU)$ of groupoids in a category $\cU$ which has finite limits. Our basic example for $\cU$ is the category $\Top$ of topological spaces. 
\begin{lem}\label{ioiowed}
The category $\groupoids(\cU)$ admits  standard models of all $2$-categorical fiber products.
\end{lem}
\proof
The objects and morphisms of the standard model of a fiber product in $\groupoids(\cU)$ can be expressed in terms of fiber products in $\cU$.\hB

\begin{lem}\label{ll1}
In $\groupoids(\cU)$ equalizers exist for any pair of maps.
\end{lem}
\proof
We observe that  $\groupoids(\cU)$ has a final object and admits $2$-categorical fiber products (Lemma \ref{ioiowed}).
In fact, the limit of the empty diagram in $\cU$ is the final object $*$ of $\cU$. The groupoid $*\Rightarrow *$ is the final object in $\groupoids(\cU)$.  \hB 

\subsubsection{}

Let $\cC$ be as in \ref{ifext}. We consider a diagram (\ref{iiuwedwnn}).
\begin{lem}\label{zuuzedzuwdwed}
We have a natural isomorphism $I(A\times_CB)\cong IA\times_{IC}IB$,
where we use standard models for the fiber products.
\end{lem}
\proof
We only have to check this for $\cC=\Cat$ since everything
can be stated in terms of Hom-categories. We let $\tilde \cD$
be the category freely generated by two objects $0$, $1$,
and two isomorphisms from $0$ to $1$, see \ref{finset-cats}.
Then we  have an isomorphism $IA \cong \Hom(\tilde \cD,A)$, see
also Lemma \ref{ixgahs} in the case of groupoids.
Since standard fiber products commute with the cotensor structure the claim follows.

%We consider the category  
%$D$ as in (\ref{kkwkehdwe}). Since both the inertia and the fibre %products are limits over $D$-diagrams in $\cC$ we can from a $D\times %D$-diagram in $\cC$ such that the left- and right-hand sides of the %asserted equivalence are obtained by taking the $2$-limits over the %factors one after the other, in the two possible orders. The asserted %equivalence is a consequence of  Proposition \ref{hjhjhjhqa} which says %that both iterated $2$-limits are equivalent to the $2$-limit over %$D\times D$.\hB 

\subsection{Loops}\label{ttzwtzetzwiu}

In a $2$-category of groupoids $\gpd(\cU)$  or stacks $\St(\cS)$ the preferred model (see \ref{zuzuuzewd}) of the inertia $IX$ (see Definition \ref{iooiwed}) of $X$ is quite complicated. The goal of the present Subsection is the construction of a simpler model of $IX$ which we call the loop object $LX$.

\subsubsection{}

We start with the case of $\gpd(\cU)$.
Let us assume that $\cU$ is tensored and cotensored over $\Sets$.
The cotensor functor will be denoted by 
$$\uHom:\Sets^{op}\times \cU\to \cU\ .$$
Using the existence of finite limits in $\cU$ we extend this functor to a bifunctor
$$\uHom_{\Cat}:(\Sets^{fin}-\Cat)\times (\cU-\Cat)\to (\cU-\Cat),$$
where for a category $\cA$ with finite limits we write $(\cA-\Cat)$ for the
$2$-category of category objects in $\cA$, and $\Sets^{fin}$ is the category of finite sets.

\subsubsection{}

Let $X\in \groupoids(\cU)\subset (\cU-\Cat)$ be a groupoid in $\cU$.
We consider the category
$$\cD:=\xymatrix{\bullet_0\ar@/^1pc/[r]^\alpha\ar@/_1pc/[r]_\beta&\bullet_1}\in (\Sets^{fin}-\Cat)\ .$$ Since $X$ is a groupoid,
$\uHom_{\Cat}(\cD,X)\in (\cU-\Cat)$ is again a groupoid in $\cU$.

\begin{lem}\label{ixgahs}
We have a natural isomorphism
$$IX\cong \uHom_{\Cat}(\cD,X)\ .$$
\end{lem}
\proof
We insert the standard model of the $2$-categorical fibre product of $\groupoids(\cU)$ into the definition of the equalizer in the special case that $f=g=\id_X$. 
Then the assertion becomes obvious.
\hB

\subsubsection{} \label{finset-cats}

Later we will have the freedom to replace groupoids by equivalent groupoids.
We let $\tilde \cD$ be the category obtained from $\cD$ by adjoining
inverses. Since $X$ is a groupoid we have
$$\uHom_{\Cat}(\cD,X)\cong \uHom_{\Cat}(\tilde \cD,X)\ .$$
We now consider the category $\cL$ with one object $*$ and infinite cyclic automorphism 
group generated by $\sigma$
$$\xymatrix{{*}\ar@(ul,dl)[]|{\sigma}}\ .$$
Then we have a natural functor
$i:\cL\to \tilde \cD$ which maps 
$*$ to $\bullet_0$ and $\sigma$ to $\beta^{-1}\circ \alpha$.
This is an equivalence of categories. It induces an equivalence of groupoids 
$$\uHom_{\Cat}(\cD,X)\cong \uHom_{\Cat}(\tilde \cD,X)\stackrel{i^*}{\to} \uHom_{\Cat}(\cL,X)\
.$$

\begin{ddd}\label{uhdwc}
The groupoid
$LX:=\uHom_{\Cat}(\cL,X)$ will be called the loop groupoid of $X$.
\end{ddd}
Note that we have an equivalence of groupoids
\begin{equation}\label{eqaass}
IX\to LX\ .
\end{equation}
If $f:X\to Y$ is a morphism in $\groupoids(\cU)$,
then composition with $f$ functorially induces a morphism
$Lf:LX\to LY$.

\subsubsection{}

It is easy to describe the objects and morphisms of the loop groupoid $LX$
explicitly.
\begin{lem}
The objects $LX^0$ and morphisms $LX^1$ of $LX$ are given by the following fibre products in $\cU$.
\begin{equation}\label{enull1}
\xymatrix{LX^0\ar[r]\ar[d]^\delta &X^1\ar[d]^{(s,r)}\\ 
X^0\ar[r]^{diag}&  X^0\times X^0}
\end{equation}
\begin{equation}\label{enull2}
\xymatrix{LX^1\ar[d]^s\ar[r]&X^1\ar[d]^s\\
LX^0\ar[r]^\delta&X^0}
\end{equation}
The range map is given (in the language of elements) by the map
$$r((x,\gamma),\mu):=(r(\mu),\mu\circ \gamma\circ \mu^{-1})\ .$$
\end{lem}

We will give another description of $LX^1$ which turns out to be useful later.
We define $P$ by the cartesian diagram
\begin{equation}\label{caer}
\xymatrix{ P\ar[d]^{(p,q)}\ar[r]^k&X^1\ar[d]^{s,r}\\
LX^0\times LX^0\ar[r]^{\delta,\delta}&X^0\times X^0 }\ ,
\end{equation}
The composition of $X$ induces a map $m:P\to X^1$ defined in the language of objects by
$$((x_0,\gamma_0),(x_1,\gamma_1),\mu)\mapsto \gamma_1^{-1} \circ \mu
\circ \gamma_0\circ \mu^{-1}\ .$$
\begin{lem}\label{wuiiuiuwed}
We have a cartesian diagram
$$\xymatrix{LX^1\ar[d]^j\ar[r]^i&X^0\ar[d]^{1}\\P\ar[r]^m&X^1}\ ,$$
where $j:=(s,r)$ and $i:=\delta\circ s$. 
\end{lem}
\proof
Consider an object $T\in \cU$. A map
$f:T\to LX^1$ is uniquely determined by a pair $(u,v)$, $u:T\to LX^0$ and $v:T\to X^1$ such that
$\delta\circ u=s\circ v:T\to X^0$. The map $u$ is given by pair $(a,b)$ of maps with $a:T\to X^0$ and $b:T\to X^1$ such that $s\circ b=r\circ b=a$.   We see that $u$ is completely determined by $b$. Note that $\delta\circ u=s\circ b=s\circ v$.
We have  $j\circ f=((s\circ b,b), (r\circ v, v\circ b\circ v^{-1}),v)$  and observe that $m\circ j\circ f=1\circ i\circ  f$.
This  construction is natural in $T\to LX^0$ and therefore determines a map
$LX^1\to P\times_{X^1} X^0$.

Consider now a map $g:T\to LX^1\to P\times_{X^1} X^0$ given by a pair
$(x,y)$ of maps $x:T\to P$ and $y:T\to X^0$ such that $m\circ x=1\circ y$.
The pair
$(p\circ x,k\circ x)$ satisfies $\delta\circ p\circ x=s\circ k\circ x$ and therefore defines a map
$f:T\to LX^1$. Again, the construction is functorial in $g$ and defines  a map
$P\times_{X^1} X^0\to LX^1$. 

We leave it to the reader to check that these maps are inverses to each other.
\hB

\subsubsection{}\label{uiwrewd}

Let $X\in \groupoids(\cU)$  and $LX$ be its loop groupoid.
Evaluation at the unique object $*$ of $\cL$ induces a functor
$LX\to X$. Therefore $LX$  can naturally be considered 
as an object of $\groupoids(\cU)/X$ (see \ref{uiduwed}). Note that a morphism in this category is a diagram
$$\xymatrix{Y\ar[rr]\ar[dr]&\ar@{:>}[d]&Z\ar[dl]\\&X&}\ ,$$
and a $2$-morphism between two such maps is a $2$-morphism
$f \leadsto g$ between the given $1$-morphisms $f,g: Y \to Z$ commuting
with the $2$-morphisms.

We will now consider group objects in $\groupoids(\cU)/X$.
They together with their products (i.e. fiber
products over $X$) will lie in a subcategory which
is equivalent as a $2$-category to a $1$-category,
so it will not be a problem to formulate what
we mean by a group object in this case.

\begin{lem}\label{wreuzhjsadc}
The loop groupoid $LX$ has a natural structure of a group object in $\groupoids(\cU)/X$.
\end{lem}
\proof
We consider the category
$\cE\in (\Sets-\Cat)$ pictured by
$$\xymatrix{\bullet_0\ar@(ul,dl)[]|a\ar[r]^b&\bullet_1 \ar@(dr,ur)[]|c }\ ,$$
where $a,c$ generate inifnite semigroups.
By $\tilde \cE$ we denote the category obtained from $\cE$ by adjoining inverses.
Then we observe that in the $2$-category $\groupoids(\cU)$
$$LX\times_X LX\cong \uHom_{\Cat}(\cE,X)\cong \uHom_{\Cat}(\tilde \cE,X) \ .$$
We define a functor
$j:\cL\to \tilde \cE$ which maps
$*$ to $\bullet_0$ and $\sigma$ to $b^{-1}\circ c\circ b\circ a$.
The pull-back $$LX\times_XLX\cong \uHom_{\Cat}(\tilde \cE,X)\stackrel{j^*}{\to}\uHom_{\Cat}(\cL,X)\cong LX$$
induces the composition law.  We leave it to the reader to write out the inverse, the unit and the remaining necessary verifications.
\hB

\newcommand{\Stacks}{{\tt St}}
\subsubsection{}\label{equforstac}

Let $\bS$ be a Grothendieck site. Then we can consider the category of presheaves of sets $\Pr\bS$. It is closed under taking arbitrary small limits. The $2$-category of strict prestacks $\PSt^{strict}\bS$ on $\bS$ is by definition the category
$\groupoids(\Pr\bS)$. By Lemma \ref{ll1} in $\PSt^{strict}\bS$ equalizers exist for all pairs of maps.

The catgeory  $\Pr\bS$ is tensored and cotensored over $\Sets$. Hence we can apply the construction of the loop groupoid in $\PSt^{strict}\bS$.
 We now consider the full $2$-subcategory of strict stacks $\St^{strict}\bS\subset \PSt^{strict}\bS$ of stacks on $\bS$. Recall that a stack is a prestack which satisfies descend conditions for objects and morphisms. This subcategory is closed with respect to
 $2$-limits and preserved by the cotensor structure. 
 For all pairs of maps  in the category  $\St^{strict}\bS$ the equalizer exists by  Lemma \ref{ll1}. Moreover, the loop object of a stack is again a stack.

While a strict prestack is a strict $2$-functor $\bS^{op}\to \groupoids(\Sets)$, a prestack is a (in general non-strict)
$2$-functor $\bS^{op}\to \groupoids(\Sets)$, i.e it preserves compositions of morphisms in $\bS$ up $2$-morphisms which satisfy coherence conditions for triple compositions as indicated
in \ref{2-functors}.
 The category of stacks is again a full subcategory of the category of prestacks on $\bS$ which satisfy certain descend conditions. Note that
$\PSt\bS$ is cotensored over $(\Sets-\Cat)$, i.e. we have a bifunctor
$$\uHom_{\Cat}:(\Sets-\Cat)\times \PSt\bS\to \PSt\bS\ .$$
This structure is induced by the corresponding cotensor structure of $(\Sets-\Cat)$, i.e.
for a category $\cD\in \Sets-\Cat$ and a prestack $X$ the value of  $\uHom_{\Cat}(\cD,X)$ on $U\in \bS$ is given 
by $$\uHom_{\Cat}(\cD,X)(U):=\uHom_{\Cat}(\cD,X(U))\ , $$ where $X(U)\in \Sets-\Cat$.
If $X$ is a stack, then  $\uHom_{\Cat}(\cD,X)$ is also a stack.

The $2$-categorical fibre product of (pre)stacks is given objectwise in $\bS$ by the $2$-categorical fibre-product in
$\groupoids(\Sets)$. Therefore, Lemma \ref{ixgahs} remains true in the categories $\PSt\bS$ and $\St\bS$. We can furthermore define the loop (pre)stack
$LX$ of a (pre)stack as in Definition \ref{uhdwc} and  (\ref{eqaass})  still induces an equivalence of (pre)stacks
$$IX\to LX\ .$$ Finally, Lemma \ref{wreuzhjsadc} holds in the sense, that for a (pre)stack $X$ the loops
$LX$ form a group object in the category of (pre)stacks over $X$.

\subsubsection{}

Like Lemma \ref{zuuzedzuwdwed} in the case of inertia stacks
we have
\begin{lem}\label{zuuzuzwed}
The inertia  functor preserves standard $2$-cartesian diagrams.
\end{lem}

\subsection{Loops of topological stacks}\label{wsdwkjdwds}

\subsubsection{}

We consider the small site $\Top$ of topological spaces and open coverings.
Let $\Stacks\Top$ be the $2$-category of stacks in topological spaces. By the observations \ref{equforstac} we can form the loop stack $LX$ of a stack
$X\in \Stacks\Top$. 
In the present subsection we show that taking loops preserves topological stacks. Furthermore we show that taking loops commutes with the classifying stack functor from topological groupoids to stacks in topological spaces. We use the latter result in order to verify that $LX$ for an orbispace is what is called the orbispace of twisted sectors in the literature.

\subsubsection{}

We refer to \cite{heinloth}, \cite{math.AG/0503247} and also to \cite{math.GT/0508550} for details about stacks (in topological spaces). Topological spaces are considered as stacks via the Yoneda embedding.
A map $a:A\to X$ from a topological space to a stack $X$ is called an atlas if it is representable, surjective and admits local sections. 
A topological stack is a stack which admits an atlas. 
We shall show that taking loops preserves topological stacks.
\begin{lem}\label{ll3} If $X\in \Stacks\Top$ is a topological stack, then
 $LX$ is a topological stack.
\end{lem}
Let  $a:A\to X$ be an atlas of $X$. Then we define a space $W$ by the pull-back diagram  
$$\xymatrix{W\ar[r]\ar[d]^w&A\times_XA\ar[d]^{(\pr_1,\pr_2)}\\
A\ar[r]^{\diag}&A\times A}\ .$$
%The map $\sigma$ induces a $2$-isomorphism $\psi:a\circ w\leadsto a\circ w$.
% It is easy to see that $(w:W\to A,\psi)$ represents the equalizer of the pair
% $a,a:A\rightarrow X$. 
% We have a diagram
% \begin{equation}\label{e1ddkl}
% \xymatrix{&W\ar[dl]_{a\circ w}\ar[dr]^{a\circ w}&\\X\ar@{:>}[rr]^{(\id_X)_*\psi}\ar[dr]^{\id_X
% }&&X\ar[dl]_{\id_X}\\&X&}\in \cE(\id_X,\id_X)\ .
% \end{equation} 
We will construct a canonical map $c:W\rightarrow LX$ and show that it is an atlas of $W$.
 
The map $c:W\to LX$ is defined as follows. Let $T$ be a topological space and $(f:T\to W)\in W(T)$.
By the definition of $W$ this map is given by a pair $(g,h)$ of maps $g:T\to A$ and $h:T\to A\times_XA$ such that $\diag\circ g=(\pr_1\circ h,\pr_2\circ h)$. The map $h:T\to A\times_XA$
is given by a pair $h_1,h_2:T\to A$ and a $2$-isomorphism $\sigma:a\circ h_1\leadsto a\circ h_2$.  Combining these two facts we see that $f$ is given by a pair $(g,\sigma)$ of a map
$g:T\to A$ and a $2$-automorphism $\sigma:a\circ g\leadsto a\circ g$.
Recall that an object of $LX(T)$ is a pair $(u,\phi)$ of an object $u\in X(T)$ and an automorphism $\phi\in \Aut(u)$.
We define $c(f)\in LX(T)$ to be the object $(a\circ g,\sigma)\in LX(T)$.

We now construct a $2$-commutative diagram
\begin{equation}\label{enull3}
\xymatrix{W\ar[d]^w\ar[r]^c&LX\ar[d]^i\\
A\ar@{:>}[ur]^\phi\ar[r]^a&X}
\end{equation}
by defining $\phi$ is follows. As above let $(f:T\to W)\in W(T)$ be given by a pair $(g,\sigma)$.
In $X(T)$ we have the equalities $i\circ c(f)=i(a\circ g,\sigma)=a\circ g$ and $a\circ w(f)=a\circ g$.
Therefore we can define $\phi(f):=\sigma$.

We claim that the diagram (\ref{enull3}) is $2$-cartesian. In order to see this let as above $T$ be a space and consider a triple $(u,v,\theta)$ consisting of maps
$u:T\to A$, $v:T\to LX$ and  a $2$-isomorphism $\theta:a\circ u\leadsto i\circ v$.
To this data we must associate a unique pair of maps $(f,\psi)$ of a map $f:T\to W$ and a $2$-isomorphism $\psi:c\circ f\leadsto v$ such that 
$$\xymatrix{T\ar[ddddr]^u\ar[rdd]^f\ar[rrrdd]^v&&&&&\\&&&&\\&W\ar@{:>}[ur]^{\psi}\ar[dd]^w\ar[rr]^c&&LX\ar[dd]^i\\&&&&\\&A\ar@{:>}[uurr]_\theta^\phi\ar[rr]^a&&X}$$
commutes.
The map $v$ is given by a pair  $(i\circ  v,\kappa)$ of an object $i\circ  v\in X(T)$ and an automorphism  $\kappa\in \Aut(i\circ  v)$. Using the description of maps $T\to W$ obtained above we define $f:T\to W$ as the map which corresponds to  the pair $(u,\theta^{-1}\circ \kappa\circ \theta)$ of an object
$u:T\to A$ and the automorphism  $\theta^{-1}\circ \kappa\circ \theta:a\circ u\leadsto a\circ u$. We furthermore define
 $2$-isomorphism 
$\psi:c\circ f=(a\circ u,\theta^{-1}\circ \kappa\circ \theta)\stackrel{\kappa^{-1}\circ \theta}{\leadsto}
 (i\circ v,\kappa)=v$. Observe that $\psi$ is uniquely determined by the condition that
$i(\psi)\circ \phi(f)=\theta :a\circ w\circ f=a\circ u\to i\circ v$. This equality  indeed holds for our construction since  $\phi(f)=\theta^{-1}\circ \kappa\circ \theta$ and $i(\psi)=\kappa^{-1}\circ \theta$.  This finishes the proof of the claim.

Since $A\rightarrow X$ is an atlas the map $a$ is representable, surjective and admits local sections. These properties are preserved under pull-back. It follows that
$c:W\rightarrow LX$ is representable, surjective and admits local sections, too. Therefore it is an atlas of $LX$.
\hB

% \subsubsection{}
% 
% \begin{lem}
% If $f:X\to Y$ is a representable map between topological stacks, then the induced map
% $Lf:LX\to LY$ is representable.
% \end{lem}
% \proof

\subsubsection{}

A topological groupoid $\cG$ is a groupoid object in $\Top$. It  represents the stack of $\cG$-principal bundles $\cB\cG$. 
If $A\rightarrow X$ is an atlas of a topological stack, then we form the topological groupoid $\cA:A\times_XA\Rightarrow A$. The stack of $\cA$-principal bundles is equivalent to $X$. We can define an equivalence $X\rightarrow \cB\cA$ which maps
$(T\to X)\in X(T)$ to  $(T\times_XA\rightarrow T)\in \cB\cA$ (we omit to write the action of $\cA$ on that space over $T$).

\subsubsection{}

Observe that finite limits in $\Top$ exist, and that $\Top$ is tensored and cotensored over $\Sets$. Therefore by \ref{ll1} for any pair of maps in $\groupoids(\Top)$ an equalizer exists.  Furthermore, we can form the loop groupoid $L\cA$ of a topological groupoid $\cA$.

\subsubsection{}

Let $A\rightarrow X$ be the atlas of a topological stack, and let $\cA\in \groupoids(\Top)$ denote the associated topological groupoid. 

\begin{lem}\label{ll4}
We have a natural equivalence of stacks
$LX\cong \cB L\cA$.
\end{lem}
\proof
Let $W\to LX$ be as in the proof of Lemma \ref{ll3}. Then can form $\cW:W\times_{LX}W\Rightarrow W$. If we show that $\cW\cong L\cA$, then the assertion follows.

>From (\ref{enull1}) we get  $W\cong (L\cA)^0$.
Next we calculate  using (\ref{enull2})
\begin{eqnarray*}
W\times_{LX}\times W&\cong&(A\times_XLX)\times_{LX}(A\times_XLX)\\
&\cong&LX\times_X(A\times_XA)\\
&\cong&(LX\times_XA)\times_A(A\times_XA)\\
&\cong&(L\cA^0)\times_A \cA^1\\
&\cong&(L\cA)^1\ .
\end{eqnarray*}
These isomorphisms are compatible with the groupoid structures.
\hB 
 
\subsubsection{}

The following result was also shown in \cite[Cor. 7.6]{math.AG/0503247}.
\begin{lem}\label{rew}
If $X$ is a topological stack, then $LX\to X$ is representable.
\end{lem}
\proof
We must show that for all spaces $T$ and maps $T\to X$ the fibre product $T\times_XLX$ is equivalent to a space. It suffices to verify this in the case that $T$ is an atlas.

We choose an atlas $A\to X$. The assertion then follows from the following two facts:
\begin{enumerate}
\item
 The diagram (\ref{enull3}) is cartesian.
\item   $W$ is a space.
\end{enumerate}
\hB

\subsubsection{}\label{uiiuduiwed}
Let us recall some notions related to orbispaces. Orbispaces as particular kind of topological stacks 
have previously been introduced in \cite[Sec. 2.1]{math.GT/0508550} and \cite[Sec. 19.3]{math.AG/0503247}). In the present paper  we use the set-up of   \cite{math.GT/0508550} 
but add the additional condition that an orbifold atlas should be separated. This condition is needed in order to show that the loop
stack of an orbifold is again an orbifold.

\begin{enumerate}
\item A topological groupoid $A:A^1\Rightarrow A^0$ is called separated if the identity $\eins_A:A^0\to A^1$ of the groupoid is a closed map.
\item A topological groupoid $A^1\Rightarrow A^0$ is called proper 
if $(s,r):A^1\to A^0\times A^0$ is a proper map. 
\item A topological groupoid is called {\'e}tale if
the source and range maps $s,r:A^1\to A^0$ are  {\'e}tale.
\item A proper \'etale  topological groupiod $A^1\Rightarrow  A^0$ is called very proper if 
there exists a continuous function $\chi : A^0\to [0,1]$ such that
\begin{enumerate}
\item $r:\supp(s^*\chi)\to A^0$ is proper
\item  $\sum_{y\in A^x} \chi(s(y))=1$ for all $x\in A^0$.
\end{enumerate}
\item A topological stack is called (very) proper ({\'e}tale, separated, resp.),  if it admits an atlas $A\to X$ such that the topological groupoid
$A\times_XA\Rightarrow A$ is (very) proper  ({\'e}tale, separated, resp).
 \item
An orbispace $X$ is a very proper {\'e}tale separated topological stack.
\item An orbispace atlas of a topological stack $X$ is an atlas $A\to X$ such that $A\times_XA\Rightarrow A$ is a very proper  {\'e}tale and separated  groupoid. 
\item If $X,Y$ are orbispaces, then a morphism of orbispaces $X\to Y$
is a representable morphism of stacks.
 \end{enumerate}

\subsubsection{}

The following Lemmas illustrates the meaning of the separatedness and very properness condition.
\begin{lem}
Let $A:A^1\Rightarrow A^0$ be a proper \'etale groupoid. If $A^1,A^0$ are locally compact, then
$A$ is very proper.
\end{lem}
\proof
The existence of the cut-off function was shown in \cite[Prop. 6.11]{MR1671260}.

\begin{lem}\label{uduiqwdqw444}
Let $A:A^1\Rightarrow A^0$ be a topological groupoid. If $A^0$ and $A^1$ are Hausdorff spaces, then $A$ is separated. 
\end{lem}
\proof
We define the Hausdorff space $Q$ as the pull-back
$$\xymatrix{Q\ar[r]^j\ar[d]&A^1\ar[d]^{(r,s)}\\A^0\ar[r]^{\diag}&A^0\times A^0}\ .$$
The property of a map between topological spaces being a closed is preserved under pull-back.
Since $A^0$ is Hausdorff the diagonal  $\diag: A^0\to A^0\times A^0$ is a closed map.
It follows that $j:Q\to  A^1$ is a closed map.
The composition $\circ$ in $A$ gives the squaring map
$$sq:Q\stackrel{\diag}{\to}Q\times_{ A^0} Q\stackrel{\circ}{\to} Q\ .$$
Then we have a pull-back
$$\xymatrix{I\ar[d]\ar[r]^k&Q\ar[d]^{(\id_Q,sq)}\\Q\ar[r]^{\diag}&Q\times Q}\ .$$
Since $Q$ is Hausdorff, it follows that $\diag $ and hence $k$ are a closed maps.
The composition $j\circ k:I\to A^1$ of closed maps is again a closed. In a group the identity
is the unique solution of the equation $x^2=x$. It follows that $j\circ k(I)=\eins_A(A^0)$.
Therefore $\eins_A(A^0)\subseteq A^1$ is closed. 

This implies that $\eins_A:A^0\to A^1$ is a closed map.
If fact, if $K\subseteq A^0$ is a closed subset, then we define the Hausdorff space
$A^1_K\subseteq A^1$ as the pull-back
$$\xymatrix{A^1_K\ar[d]\ar[r]^v&A^1\ar[d]^{(r,s)}\\K\times K\ar[r]^u&A^0\times A^0}\ .$$
Since $u$ (the obvious embedding) is a closed map, so is $v$. We apply the discussion above to  the restricted groupoid   $A^1_K\Rightarrow K$ with identity $\eins_{A_K}:K\to A^1_K$ in order to show that $\eins_{A_K}(K)\subseteq A^1_K$ is closed. Hence $\eins_A(K)=v(\eins_{A_K}(K)))\subseteq A^1$ is closed.  
\hB

\subsubsection{}

\begin{lem}\label{ll56}
If $X$ is an orbispace, then $LX$ is an orbispace and $LX\rightarrow X$ is a morphism of orbispaces.
\end{lem}
\proof
We choose an orbispace atlas $A\to X$. The associated groupoid $\cA:A\times_XA\rightarrow A$ is {\'e}tale, proper. and separated. In order to show that $LX$ is an orbispace it suffices to show by Lemma \ref{ll4} that $L\cA$ is {\'e}tale, proper and separated.

The property of a map between topological spaces being \'etale is preserved under pull-back. By  (\ref{enull2})  the fact that $s:\cA^1\rightarrow \cA^0$ is {\'e}tale therefore implies that $s:(L\cA)^1\rightarrow (L\cA)^0$ is {\'e}tale. Using the inversion homeomorphism
$I:(L\cA)^1\rightarrow (L\cA)^1$ we can express the range map in terms of the source map: $r=s\circ I$. This implies that $r:(L\cA)^1
\to (L\cA)^0$ is \'etale, too. We thus have shown that $L\cA$ is {\'e}tale.

We consider the pull-back
$$\xymatrix{P\ar[r]\ar[d]^j&\cA^1\ar[d]^{(r,s)}\\(L\cA)^0\times (L\cA)^0\ar[r]&\cA^0\times \cA^0}$$
(compare (\ref{caer}) ). The property of a map between topological spaces being proper is also preserved by pull-backs.
Therefore $j:P\to (L\cA)^0\times (L\cA)^0$ is a proper map.
The image of $\eins_{\cA}:\cA^0\to \cA^1$ is closed. By Lemma \ref{wuiiuiuwed} we can write $(L\cA)^1$ as a closed subspace $(L\cA)^1:=m^{-1}(\eins_{\cA}(\cA^0))\subset P$. In general, the restriction of a proper map to a closed subspace is still proper. Since the restriction of  $j$ to the  closed subspace $(L\cA)^1\subset P$ is exactly
$(r,s):(L\cA)^1\to (L\cA)^0\times(L\cA)^0$ we see that the groupoid $\cA$ is proper.\footnote{It is because of this argument that in addition to the conditions used in \cite{math.GT/0508550} we  require an orbispace atlas to be separated.}

We now show that $L\cA$ is very proper. Since $\cA$ is very proper there exists a continuous function
$\chi:\cA^0\to [0,1]$ such that $r:\supp(s^*\chi)\to \cA^0$ is proper and
$\sum_{y\in \cA^x} \chi(s(y))=1$ for all $x\in \cA^0$. Let $i:L\cA\to \cA$ be the canonical map.
Then $i^*\chi:L\cA^0\to [0,1]$ has corresponding  properties for the groupoid $L\cA$.

Finally we show that $\eins_{L\cA}:(L\cA)^0\to (L\cA)^1$ is a closed map.
By definition we have the cartesian square
$$\xymatrix{(L\cA)^1\ar[r]\ar[d]&\cA	^1\ar[d]\\(L\cA)^0\ar[r]&\cA^0}\ .$$
Therefore we have an embedding as a subspace
$(L\cA^1)\subset (L\cA)^0\times \cA^1$. Let $K\subseteq (L\cA)^0$ be a closed subset. 
Then we can write
$\eins_{L\cA}(K)=(L\cA^1)\cap (K\times i(\cA^0))$. Since $\cA$ is separated the subspace
$ (K\times \eins_{\cA}(\cA^0))\subseteq  (L\cA)^0\times \cA^1$ is closed.
Therefore $\eins_{L\cA}(K)\subset (L\cA)^1$ is closed, too.

In order to be a map of orbispaces $LX\to X$ must be representable. This is Lemma \ref{rew}.
\hB

\subsubsection{}\label{smoothworld}

We can replace the site of topological spaces $\Top$ by the site of smooth manifolds $\Mf^\infty$. We will call the corresponding stacks stacks in smooth manifolds. A map $A\to X$ from a manifold 
to a stack in smooth manifolds is called an atlas if it is representable, surjective and smooth (i.e. submersion). A stack in smooth manifolds which admits an atlas is called a smooth (or differentiable) stack. An orbifold is a proper {\'e}tale smooth stack in smooth manifolds.
Since manifolds are Hausdorff a smooth stack is

% \begin{lem}
% If $X$ is a smooth stack, then $LX$ is a smooth stack. 
% If $X$ is an orbifold, then so is $LX$.
% \end{lem}
% \proof
The obvious problem to extend the proof of Lemma \ref{ll3} from topological spaces to smooth manifolds is that in smooth manifolds fibre products only exist under appropriate transversality conditions. In fact, the map $(\pr_1,\pr_2):A\times_X A\to A\times A$ is in general  not transverse to the diagonal $\diag:A\to A\times A$.

But it is still true that the loop stack of an orbifold is an orbifold. Proofs of this fact can be found e.g. 
in  \cite{MR0474432}, \cite{MR1993337},  \cite{MR2104605}.
Note that for smooth stacks $LX\to X$ is in general neither smooth nor representable.

\subsection{Loops and principal bundles}\label{hdjjhsjcoo}

\subsubsection{}

Let $G$ be a topological group. The classifying stack $\cB G$ of $G$-principal bundles is given as a quotient stack
$\cB G:=[*/G]$ of the action of $G$ on the one point space ${*}$ \cite[Example 2.5]{heinloth}. The map $*\to \cB G$ is an atlas and we have a canonical cartesian diagram
$$\xymatrix{G\ar[d]\ar[r]&{*}\ar[d]\\{*}\ar@{:>}[ur]\ar[r]&\cB G}\ .$$
Hence this atlas gives rise to the groupoid
$\cG:G\Rightarrow *$. We see that $L\cG$ is the groupoid $G\times G\Rightarrow G$
of the action of $G$ on itself by conjugations. Therefore by Lemma \ref{ll4} we have $L\cB G\cong [G/G]$.

\subsubsection{}

A $G$-principal bundle over a space $Y$ is by definition an object of $p\in \cB G(Y)$, or equivalently, by Yoneda's Lemma, a map
$p:Y\to \cB G$. The underlying map of spaces $P\to Y$ fits into the cartesian diagram
$$\xymatrix{P\ar[r]\ar[d]&{*}\ar[d]\\Y\ar@{:>}[ur]\ar[r]^p&\cB G}\ .$$
We adopt the same definition for a $G$-principal bundle over a stack $Y$. In this case the underlying map $P\to Y$ is a representable map. 

\subsubsection{}\label{cohjhjhjsd}

Let $a:A\to Y$ be an atlas such that the pull-back of the principal bundle $p:Y\to \cB G$ admits a trivialization. A trivialization is a lift $t$ in the diagram $$\xymatrix{&&{*}\ar[d]\\A\ar@{.>}[urr]^t\ar[r]^a&Y\ar@{:>}[ur]\ar[r]^p&\cB G}\ .$$
The cocycle associated to the atlas $a$ and the trivialization
$t$ is the induced map
$$\Phi_{a,t}:A\times_YA\to *\times_{\cB G} *\cong G\ .$$
Let $\cA:A\times_YA\Rightarrow A$ be the groupoid determined by the atlas and $\cA^\bullet$ denote the associated simplicial space. Let 
$$C^\bullet(\cA;G):=C(\cA^\bullet,G)\ ,\quad \delta:C^\bullet(\cA;G)\to C^{\bullet+1}(\cA;G)$$
be the associated cochain complex (the part in degree $> 2$ is only defined if $G$ is abelian).
Then $\Phi_{a,t}\in C^1(\cA,G)$ is closed, i.e. it satisfies $\delta\Phi_{a,t}=0$.
We refer to \cite[Sec.2]{heinloth} for a  description of $G$-principal principal bundles in terms of cocycles.

\subsubsection{}

Let $p:Y\to \cB G$ be a $G$-principal bundle over a stack $Y$. We apply the loop functor and get the map
$Lp:LY\to L\cB G\cong [G/G]$. It is a homomorphism over the map $Y\to \cB G$.
If $G$ is abelian, then it induces a homomorphism 
\begin{equation}\label{kjjkkjsdsdssds}
h:LY\to G\ .
\end{equation}

% We start this subsection with an explanation of the notion of a $G$-principal bundle in topological stacks. At the moment let $G$ be a space. A representable map $f:X\to Y$ between topological stacks is a locally trivial fibre bundle in topological stacks with fibre $G$ if
% there exists an atlas $A\to Y$ and a diagram
% \begin{equation}\label{jkkjkjeded}\xymatrix{A\times G\ar[rd]^{\pr_A}\ar[r]^\cong&B\ar[d]\ar[r]&X\ar[d]\\&
% A\ar@{:>}[ur]\ar[r]&Y}\ ,\end{equation}
% where the square is cartesian.
% We further get a natural isomorphism
% $$\Phi:A\times G\cong\pr_0^*B\stackrel{\sim}{\to} \pr_1^*B\cong A\times G$$
% over $A\times_YA$ which satsifies a cocycle condition over $A\times_YA\times_YA$.
% Here $\pr_i:A\times_YA\to A$, $i=0,1$, are the two projections on the factors.
% Vice versa,  given the cocyle $\Phi$ we can define a locally trivial $G$-bundle over $Y$ by glueing.
% Assume now that $G$ is a topological group. If $\Phi$ is an isomorphism of right $G$-spaces, then $X\to Y$ aquires the structre of a $G$-principal bundle. This is our constructive definition of a $G$-principal bundle. A characterization opposite to this constructive point of view  is discussed in \cite[Sec.2]{heinloth}. Observe that the action of $m:X\times G\to X$ induces an isomorphism
% \begin{equation}\label{uiuidednededeee}
% (\pr_X,m):X\times G\stackrel{\sim}{\to} X\times_YX\ .
% \end{equation}

\subsubsection{}

In the following we give a heuristic description of this homomorphism.
Let $f:P\to Y$ be the underlying map of stacks of the principal bundle. Furthermore let $i:LY\to Y$ denote the canonical map. 
  For a point $y\in Y$ we get an action of the group  $i^{-1}(y)$ on the fibre $f^{-1}(y)$.
If $\gamma\in i^{-1}(y)$ and $x\in f^{-1}(y)$, then
$\gamma  x= x h(\gamma)$.  
On the left-hand side, $(\gamma,x)\mapsto \gamma x$ denotes the action of $i^{-1}(y)$ on $f^{-1}(y)$. On the right-hand side $(x,g)\to xg$ is the $G$-action on $P$ given by the principal bundle structure. We see again, that the restriction $h_{|i^{-1}(y)}:i^{-1}(y)\to G$ is a homomorphism for all $y\in Y$.

\subsubsection{}\label{wuiefdhwe}
Assume that we have chosen an atlas $a:A\to Y$ and a trivialization $t$ as in \ref{cohjhjhjsd}.
Let $\cA:A\times_YA\Rightarrow A$ be the associated groupoid.
Then we get an induced map
$h_a:L\cA\rightarrow G$. It is equal to the restriction of the cocycle $\Phi_{a,t}$ to $(L\cA)^0\subseteq \cA^1$, i.e. we have the equality
\begin{equation}\label{weuiuwefuiewdw}
h_a=(\Phi_{a,t})_{|(L\cA)^0}\ .
\end{equation}
The cocycle $h_a$ is closed, i.e. $\delta h_a=0$, and it represents the function $h\in C(LY;G)$ under the identification  $H^0(L\cA;G)=C(LY,G)$.
Another interpretation of (\ref{weuiuwefuiewdw}) is as the equality  $h_a=\tr[\Phi_{a,t}]$, where $[\Phi_{a,t}]\in H^1(\cA;G)$ is the cohomology class represented by $\Phi_{a,t}$, and $\tr:C^{\bullet+1}(\cA;G)\to C^{\bullet}(L\cA;G)$ is the  transgression chain map
defined in  \cite{math.AT/0605534}, \cite{math.AT/0307114}, \cite{math.KT/0604160}.

\subsubsection{}\label{ejkjkjkked666}

Let $G$ be an abelian  topological group. In the following Lemma we
will assume that for all $n\in \nat$ the subspace of $n$-torsion points $$\Tors_n(G):=\{g\in G|g^n=1\}\subseteq G$$ is discrete. This is  a non-trivial assumption which, for example,  is not true for the  topological group $\prod_{\nat} \Z/n\Z$.   Let $G^\delta$ denote the group $G$ with the discrete topology.
Let $p:Y\to \cB G$ be a $G$-principal bundle. 
 
\begin{lem}\label{ewhjd}
 If $Y$ is an orbispace and the subsets $\Tors_n(G)\subseteq G$ are discrete for all $n\in \nat$, then the map $h:LY\to G$  (defined in (\ref{kjjkkjsdsdssds})) factors over $G^\delta$.
\end{lem}
\proof
We must show that for all spaces $T$ and maps $w:T\to LY$ the composition $h\circ w:T\to G$ is locally constant.
We choose an orbifold atlas $A\to Y$ which gives rise to a very proper separated \'etale groupoid $\cA:A\times_YA\Rightarrow  A$.
 
We consider a point $t\in T$. There exists a neighbourhood $t\in U\subseteq T$ which admits a lift
$$\xymatrix{U\ar@{.>}[r]^{\tilde w}\ar[d]&\cA^0\ar[d]\\T\ar@{:>}[ur]^\sigma\ar[r]^{w\circ i}&Y}\ .$$ 
By Lemma \ref{ll4} we have the $2$-cartesian square in the following diagram:
$$\xymatrix{U\ar@{.>}[dr]^{v}\ar@/^0.7cm/[drr]^{w}\ar@/^-0.7cm/[ddr]_{\tilde w}&&\\&(L\cA)^0\ar[d]\ar[r]&LY\ar[d]^i\\&
\cA^0\ar@{:>}[ur]^\sigma\ar[r]&Y}\ .$$
We get an induced map 
$v:U\to L\cA^0\subseteq \cA^1$ such that $\tilde w=s\circ v$. Let $a:=\tilde w(t)\in \cA^0$ so that $v(t)\in \cA_a^a$. Since the groupoid $\cA$ is proper the group $\cA_a^a$ is finite.  Hence there exists an $n\in \nat$ such that $v(t)^n=\id_a$.
The map $v^n$ fits into the diagram
$$\xymatrix{{\{t\}}\ar[r]\ar[d]&\cA^1\ar[d]^s\\U\ar[ur]^{v^n}\ar[r]^{\tilde w}&\cA^0}\ .$$
Note that the map $U\ni u\mapsto \id_{\tilde w(u)}\in \cA^1$ would fit into the same diagram in the place of $v^n$.
Since $s:\cA^1\to \cA^0$ is \'etale we can shrink $U$ further
such that $v^n(u)=\id_{\tilde w(u)}$ for all $u\in U$.
This implies that $h\circ w_{|U}:U\to G$ factors over the discrete subset $\Tors_n(G)\subseteq G$ an is therefore locally constant.
\hB

% Since $X\to Y$ admits local sections we can
% find an open covering $(T_\alpha\to T)_{\alpha\in I}$ of $T$ 
% such that for all $\alpha \in I$ there exists
% a section, i.e. a map $s_\alpha:T_\alpha\to X$ and a $2$-morphism $\sigma_\alpha$ such that the following diagram commutes.
% $$\xymatrix{T_\alpha\ar[r]^{s_\alpha}\ar[d]&X\ar[d]\\LY\ar@{:>}[ur]^{\sigma_\alpha}\ar[r]&Y}$$
% By \ref{hwsjhwsw} such a diagram naturally gives a  map
% $$h_{T_\alpha}:T_\alpha\to X\times_YX\stackrel{(\ref{uiuidednededeee})}{\cong} X\times G\stackrel{\pr_2}{\to} G\ .$$
% Using the assumption that $G$ is abelian we check that this map does not depend on the choice of $s_\alpha$ and $\sigma_\alpha$. Let $s_\alpha^\prime$, $\sigma_\alpha^\prime$ be another choice. Then we can write
% $s_\alpha^\prime:=s_\alpha u$ for some function $u:T_\alpha\to G$.
% 
% 
% 
% 
% 
% 
% 
% Therefore the collection of  maps $(h_{T_\alpha})_{\alpha\in I}$ defines a  global map $h_T:T\to G$.
% Since this construction is functorial in $T$ it determines a map
% $h:LY\to G$.
% Finally, it is easy to check that it is a group stack homomorphism.
% 
% If $X$ is an orbispace, then using an orbispace atlas we can write out this construction 
% in terms of proper \'etale groupoids and see that the map factors over $G^\delta$.
% \hB
% As a side remark, if $G$ is not abelian, then by a similar construction we can define a map
% $LY\to G/G$, where $G/G$ the space of conjugacy classes of $G$.

\subsubsection{}

Let $G$ be a topological abelian group such that $\Tors_n(G)\subset G$ is discrete for all $n\in \nat$.
Furthermore, let $p:Y\to \cB G$ be a $G$-principal bundle over an orbisspace $Y$ and
$h:LY\to G^\delta$ as in Lemma \ref{ewhjd}. Then we have a decomposition
$$LY\cong \bigsqcup_{g\in G} LY_g\ ,$$
where $LY_g:=h^{-1}(g)$ is formally defined by the $2$-cartesian square
$$\xymatrix{LY_g\ar[r]\ar[d]&[\{g\}/G]\ar[d]&\\LY\ar@{:>}[ur]\ar[r]^{Lp}&[G^\delta/G]\ar[r]^\cong& \bigsqcup_{l\in G}[\{l\}/G]}\ .$$

Let $f:X\to Y$ be the map of stacks underlying the principal bundle $p$. It fits
 into the cartesian diagram
\begin{equation}\label{iisuuhdsddddddd}\xymatrix{X\ar[d]^f\ar[r]&{*}\ar[d]\\Y\ar@{:>}[ur]\ar[r]^p&\cB G}\ .\end{equation}
\begin{lem}\label{hnsdkjc}
  The map
$Lf:LX\to LY$ factors over the $G$-principal bundle
$$LX\to LY_1\ .$$ 
\end{lem}
\proof
We apply the loop functor to the $2$-cartesian diagram (\ref{iisuuhdsddddddd}) and get the $2$-cartesian diagram  (see Lemma \ref{zuuzuzwed})
\begin{equation}\label{zewuzzuwezufwef}
\xymatrix{LX\ar[d]^{Lf}\ar[r]&L\{1\}\ar[d]\ar@{=}[r]&\{1\}\ar[d]\\LY\ar@{:>}[ur]\ar[r]^p&L\cB G\ar@{=}[r]&[G/G]}\ .
\end{equation}
It follows from the construction of $h:LY\to G$ that $h\circ Lf$ is the constant map with value $1\in G$. It remains to show that $LX\to LY_1$ is a $G$-principal bundle. To this end we refine the diagram (\ref{zewuzzuwezufwef}) to
$$\xymatrix{LX\ar[d]^{Lf}\ar[r]&\{1\}\ar[d]\\LY_1\ar[r]\ar[d]&[\{1\}/G]\ar[d]&\\LY\ar@{:>}[uur]\ar@{:>}[ur]\ar[r]^p&[G/G]}\ .$$
By definition of $LY_1$ the lower square is $2$-cartesian. Since the outer square is the $2$-cartesian square (\ref{zewuzzuwezufwef}) we conclude that
the upper square is $2$-cartesian.
\hB

\subsubsection{}\label{ejfwws}
Let $\Gamma$ be a finite group.
The exact segment
 $$\xymatrix{\dots\ar[r]&H^1(\Gamma;\R^\delta)\ar[r]\ar[d]^\cong&H^1(\Gamma;U(1)^\delta)\ar[d]^\cong\ar[r]_\partial^\cong &H^2(\Gamma;\Z)\ar[r]&H^2(\Gamma;\R^\delta)\ar[d]^\cong\ar[r]&\dots\\
 &0&\hat \Gamma&&0&}$$
of the Bockstein sequence in group cohomology associated to the sequence of coefficients
$$0\to \Z\to \R^\delta\to U(1)^\delta\to 0$$
gives rise to a natural identification
$$H^2(\Gamma;\Z)\cong \hat \Gamma\ ,$$ where $\hat \Gamma$ 
denote the group of $U(1)$-valued characters of $\Gamma$.

Let us consider the orbispace $[*/\Gamma]$. Then we have
$L[*/\Gamma]\cong [\Gamma/\Gamma]$, where $\Gamma$ acts on itself by conjugation.
A character  $\chi\in \hat \Gamma$ gives rise to a function
$$\bar \chi:L[*/\Gamma]\to U(1)^\delta\ ,\quad \gamma\mapsto \chi(\gamma)\ .$$

\subsubsection{}

There are various ways to define the integral cohomology of an orbispace $B$.
In order to be able to use results about the classification of $U(1)$-principal bundles over $B$ 
we use the definition \cite{math.GT/0508550}, where we define  $$H^*(B;\Z):=H^*(|\cA|;\Z)$$  using the classifying space 
$|\cA|$ of the groupoid $\cA$ associated to an orbifold atlas $a:A\to B$.
% If we choose a good atlas such that $\cA^\bullet$ is contractible in all degrees ({\bf union of contractibles?}), then
% we also have an isomorphism $H^*(B;\Z)\cong H^*(\cA;\Z)$ ({\bf definition?}).
Note that by this definition $H^*(\cB \Gamma;\Z)\cong H^*(\Gamma;\Z)$. 
In fact, if we choose the atlas $a:*\to \cB\Gamma$ and let $\cA$ be the associated groupoid, then $|\cA|$ is the standard model of the classifying space $B\Gamma$ of $\Gamma$.
\subsubsection{}\label{bxncmy}

Let $\chi\in H^2(B;\Z)$. 
In this paragraph we generalize the construction \ref{ejfwws} of the map
$\chi\mapsto \bar \chi$ to general orbispaces $B$. We start with describing the values of $\bar \chi:LB\to U(1)$ at the points of $LB$.
For the moment we do not claim any continuity property, but by Lemma \ref{jujzsdh} we see that it is continuous even if we equip $G$ with the discrete topology.

Consider a point $u:*\to LB$. It determines and is determined by a point
$p_u:*\stackrel{u}{\to} LB\to  B$ in $B$ and an element $\gamma_u\in \Aut(p_u)\cong *\times_B*$.
The element $\gamma_u$ generates a finite cyclic group $\Gamma_u$. We obtain an induced map
$\tilde u:[*/\Gamma_u]\to B$. We have  $L[*/\Gamma_u]\cong [\Gamma_u/\Gamma_u]$ and consider $\gamma_u\in [\Gamma_u/\Gamma_u]$
(or more formally, as a map $\gamma_u:*\to [\Gamma_u/\Gamma_u]$).
 We have an induced map
$L\tilde u:L[*/\Gamma_u]\to LB$ such that
$L\tilde u(\gamma_u)=u$.  
We can now define
$$\bar \chi(u):=\overline{\tilde u^*\chi}(\gamma_u)\ .$$

\subsubsection{}

Let $B$ be an orbispace.
By  \cite[Proposition 4.3]{math.GT/0508550})
the class $\chi\in H^2(B;\Z)$ classifies a $U(1)$-principal bundle $P_\chi\to B$.
In Lemma \ref{jujzsdh} we will express the corresponding function $h_\chi:LB\to U(1)^\delta$
(defined in (\ref{kjjkkjsdsdssds}) directly in terms of $\chi$.

% Let $B$ be an orbispace.
% The class $\chi\in H^2(B;\Z)$ classifies a $U(1)$-principal bundle $P_\chi\to B$
% (see \cite[Proposition 4.3]{math.GT/0508550}). By Lemma \ref{hnsdkjc} we get a map
% $h_\chi:LB\to U(1)^\delta$.
\begin{lem}\label{jujzsdh}
We have the equality
$h_\chi=\bar \chi$.
\end{lem}
\proof
The constructions of $h_\chi$ and $\bar \chi$  are natural under pull-back.
It therefore suffices to show this equality in the case that $B\cong [*/\Gamma]$ for a finite group $\Gamma$.
In this case we have  $P_\chi\cong [U(1)/_\chi \Gamma]$, where $\Gamma$ acts on $U(1)$ via the character $\chi$. By construction of $h_\chi$ we have
$h_\chi=\chi:[\Gamma/\Gamma]\to U(1)^\delta$. On the other hand, again by construction, we have $\bar \chi=\chi:[\Gamma/\Gamma]\to U(1)^\delta$. \hB 

\subsubsection{}

Here is another interpretation. Let $a:A \to B$ be a good orbifold atlas.
We can choose a trivialization $t$ of the pull-back of the $U(1)$-bundle to $A$
and get a cocycle $\Phi_{a,t}\in C^1(\cA;U(1))$.
The definition of an orbifold atlas is in particular made such that
$H^i(\cA;\R_{cont})=0$ for $i\ge 1$.\footnote{For a proof see \cite[Proposition 1]{math.DG/0008064} or  the corrected version \cite{math.GT/0508550}. In the orginal version an orbifold atlas was characterized by the property that it gives rise to a proper \'etale groupoid. In order to prove this vanishing of real continuous cohomology we added the assumption of being very proper.}
 Hence the boundary operator in cohomology associated to the sequence $0\to \Z\to \R_{cont}\to U(1)_{cont}\to 0$ induces an isomorphism $\partial:H^1(\cA;U(1)_{cont})\stackrel{\sim}{\to} H^2(\cA;\Z)\cong H^2(B;\Z)$.
Under this identification we have $\chi\cong \partial  [\Phi_{a,t}]$. Our construction of $\chi\mapsto \bar \chi$ is made such that
$\overline{\partial (\phi)}\cong \tr\phi\in H^0(L\cA;U(1))\cong  C(LB,U(1))$ for every class
$\phi\in H^1(\cA;U(1))$. In view of \ref{wuiefdhwe} this assertion is equivalent to Lemma \ref{jujzsdh}.

\subsection{Gerbes and local systems}\label{twto}

\subsubsection{}

We consider stacks in topological spaces $\St\Top$.
Let $H$ be an abelian topological group
and $f:G\rightarrow X$ 
be a topological gerbe with band $H$ over some topological stack $X$
We take loops and obtain 
$Lf:LG\rightarrow LX$. We further have a canoical map $\tilde i:LG\to G$, and $LG/G$ is a group in $\St\Top/G$ (see Lemma \ref{wreuzhjsadc}).
Since $i\circ Lf\cong f\circ \tilde i$ we get the dotted arrow
\begin{equation}\label{kk77qghhgsghvqhgsq}
\xymatrix{LG\ar@{.>}[dr]^\pi\ar[ddr]_{Lf}\ar[rrd]^{\tilde i}&&\\&G_L\ar[d]\ar[r]&G\ar[d]^f\\&LX\ar@{:>}[ur]\ar[r]^i&X}\ ,
\end{equation}
where the gerbe $G_L\to LX$ is defined by the $2$-cartesian square. One way to say that the gerbe $G\to X$ has band $H$ is as follows:\footnote{The definition given in \cite[Def. 5.3]{heinloth} expresses these properties using objects.} 
\begin{enumerate}
\item The map $\pi:LG\to G_L$ is the underlying map of an  $H$-principal bundle $G_L\to BH$. 
\item The sequence of (representable, see \ref{rew}) maps $\pi:LG\to G_L\to G$ is a central extension
of groups  
\begin{equation}\label{kkkqwu7}
G\times H/G\to LG/G\to G_L/G
\end{equation} in $\St\Top/G$ (the group stack structures of $G_L/G$ is induced from that of $LX/X$.
\end{enumerate}

\subsubsection{}

\begin{prop}\label{ewiudiwed}
There exists a canonical central extension
$$X\times H/X\to \tilde G/X\to LX/X$$ of groups in  $\St\Top/X$ whose pull-back along $G\to X$ is 
 isomorphic to  (\ref{kkkqwu7}). 
It depends functorially on the datum $G\to X$.
\end{prop}
\proof
We first go over to topological groupoids by choosing atlases. Then we construct the required extension in topological groupoids. Finally we pass back to stacks.\footnote{This argument is not satisfactory. It would be better to argue directly with stacks. But at the moment we do not know how to do this.} 

We choose an atlas $a:A\rightarrow X$ 
which admits a lift 
\begin{equation}\label{jkkjkwdwdwe64}
\xymatrix{&G\ar[d]^f\\A\ar@{.>}[ur]^b\ar[r]^a&X\ar@(l,ul)[]^\phi}
\end{equation} to an atlas of $G$.
We get topological groupoids 
\begin{eqnarray*}
\cX&:&\cX^1:= A\times_X A\Rightarrow \cX^0:= A\\
\cG&:&\cG^1:= A\times_GA\Rightarrow \cG^0:= A\ ,\end{eqnarray*}
  and a central $H$-extension 
$$\xymatrix{X^0\times H\ar[d]&&\\
\cG^1\ar[d]\ar@{=>}[r]&\cG^0\ar@{=}[d]\\\cX^1\ar@{=>}[r]&\cX^0}\ .$$

Using the description (\ref{enull1}) of the objects of $L\cX$ and $L\cG$ we get the pull-back of $H$-principal bundles
$$\xymatrix{(L\cG)^0\ar[d]\ar[r]&\cG^1\ar[d]\\(L\cX)^0\ar[r]&\cX^1}$$
Furthermore, by  (\ref{enull2}) we have the following description of morphisms $(L\cG)^1$ as a pull-back
\begin{equation}
\xymatrix{(L\cG)^1\ar[d]^s\ar[r]&\cG^1\ar[d]^s\\
(L\cG)^0\ar[r]^\delta&\cG^0}\ .
\end{equation}
We see that $(L\cG)^1$ has two commuting  $H$-actions, the first comes from the action on $(L\cG)^0$
(the principal bundle structure of the left lower corner in(\ref{enull2})), and the second comes from the action on $\cG^1$, the right upper corner in (\ref{enull2}).

\subsubsection{}\label{hjqwhdjqwdq}

We now define the groupoid $\cG_L$ corresponding to the stack $G_L$. The obvious definition would be as $L\cX\times_{\cX}\cG$, but we consider the simpler equivalent groupoid
 $\cG_L:(\cG_L)^1\Rightarrow (L\cX)^0$ where the morphisms are given by the cartesian diagram
\begin{equation}\label{hh32562e2398e83}
\xymatrix{(\cG_L)^1\ar[d]\ar[r]&\cG^1\ar[d]^{(r,s)}\\(L\cX)^1\ar[r]^{i\circ r,i\circ s}&\cX^0\times \cX^0}\ .
\end{equation}
We have a natural homomorphism of groupoids $L\cG\to \cG_L$ which is an $H$-principal bundle as expected.

\subsubsection{}\label{awzegbdsajmnd}

Observe that we can define a groupoid $\tilde \cG:\tilde \cG^1\rightarrow \tilde \cG^0=L\cG^0$ by taking the quotient of $\tilde \cG^1:=(L\cG)^1$ by the second $H$-action. In other words, we define $\tilde \cG^1$ by the cartesian diagram
\begin{equation}\label{deftg}\xymatrix{\tilde \cG^1\ar[d]\ar[r]&(L\cG)^0\ar[d]\\\cX^1\ar[r]&\cX^0}\ .\end{equation}
With the natural induced map $\tilde \cG\rightarrow L\cX$ is an $H$-principal bundle over $L\cX$.
We compose this map with $L\cX\to \cX$ and observe that
the groupoid structure on $\tilde \cG$ induces on $\tilde \cG\to \cX$ the structure of a group in groupoids over $\cX$. It fits into the central extension
$$\cX\times H\to \tilde \cG\to L\cX\ .$$
of groups in $\gpd(\Top)/\cX$.

The bundle $\tilde \cG\to L\cX$  fits into a cartesian diagram
$$\xymatrix{L\cG\ar[d]\ar[r]&\tilde \cG\ar[d]\\\cG_L\ar@{:>}[ur]\ar[r]&L\cX }\ .$$

\subsubsection{}

We now pass back to stacks. 
We interpret the $H$-principal bundle $\tilde \cG^0\to (L\cX)^0$ as an object
$(L\cX)^0\to \cB H$. The action $(L\cX)^1\times_{(L\cX)^0} \tilde \cG^0\to \tilde \cG^0$
gives the descend\footnote{Let $\cB:\cB^1\Rightarrow \cB^0$ be a topological groupoid with quotent stack $[\cB^1/\cB^0]$. Let $U$ be some stack. A descend datum
is a diagram
$$\xymatrix{\cB^1\ar[r]^s\ar[d]^r&\cB^0\ar[d]\\\cB^0\ar@{:>}[ur]\ar[r]&U}$$
which is compatible with the composition in $\cB$ in the obvious way.
We use the equivalence of the category  $\Hom([\cB^1/\cB^0],U)$  with the category of descend data.} datum for completing the following diagram by the dotted arrows:
$$\xymatrix{\tilde G^0\ar[r]\ar[d] &\tilde G\ar@{.>}[d]&\\(L\cX)^0\ar[dr]\ar[r]&[(L\cX)^0/(L\cX)^1]\ar@{.>}[d]\ar[r]^{\hspace{0.5cm}\cong}&LX\\&\cB H&}
$$ 

The $H$-principal bundle in groupoids $\tilde \cG\to L\cX$ thus gives rise to a $H$-principal bundle in topological stacks $LX\to \cB H$ with underlying map of stacks $\tilde G\to LX$.
In fact, it fits into the cartesian diagram $$\xymatrix{LG\ar[d]\ar[r]&\tilde G\ar[d]\\G_L\ar@{:>}[ur]\ar[r]&LX }$$
and the central extension
$$X\times H\to \tilde G\to LX$$ in
$\St\Top/X$.

In order to answer the question wether $\tilde G\to LX$ is well-defined up to canonical equivalence
we must study how it depends on the choice of the atlas $a:A\to X$ and its lift $(b,\phi)$ (see \ref{jkkjkwdwdwe64}). We must show that an automorphism of this datum induces the identity on
$\tilde \cG\to L\cX$. Now observe that the automorphism group of $(a,b, \phi)$ is the group of automorphisms of $b$ which induce the identity on $a$ (in order not to change $\phi$). By the definition of a $H$-banded gerbe it is given by $C(A,H)$. It acts trivially on $\tilde \cG\to L\cX$, indeed. 

Finally observe that the construction of $\tilde G\to LX$ depends functorially on $G\to X$. We leave the details to the reader.
\hB

\subsubsection{}

We now assume that the stack $X$ is an orbispace.
We further assume that  $\Tors_n(H)\subseteq H$ is discrete (compare \ref{ejkjkjkked666}). Let $H^\delta$ be the group $H$ equipped with the discrete topology.
\begin{lem}\label{flat}
The $H$-bundle $\phi:\tilde G\rightarrow LX$ admits a natural reduction of structure groups $\phi^\delta:\tilde G^\delta\rightarrow LX$ from $H$ to $H^\delta$.
\end{lem}
\proof
Let $T$ be a space and $*\in T$ be a distinguished point.
We consider the lifting problem
$$\xymatrix{{*}\ar[r]^\sigma\ar[d]&\tilde G\ar[d]\\T\ar[r]\ar@{.>}[ru]&LX}\ .$$
We must show that this problem has a unique solution after replacing $T$ by some neighbourhood of $*$, if necessary.

Using an orbispace atlas $A\rightarrow X$ we  translate to an equivalent  lifting problem for topological groupoids
$$\xymatrix{{*}\ar[r]\ar[d]&\tilde \cG\ar[d]&\\T\ar[r]^t\ar@{.>}[ru]^{\tilde t}&L\cA\ar[r]^i&\cA}\ .$$
Here we consider $T$ as a groupoid $T\Rightarrow T$ in the canonical way.
Let $\gamma:=t(*)\in (L\cA)^0\cong \cA^a_a$,
where $a:=i(\gamma)\in \cA^0=A$.
Since $\cA^a_a$ is a finite group there exists $n\in \nat$ such that $\gamma^n=\id_{\cA_a^a}$. We consider the embedding $\cA^0\subset \cA^1$ given by the identities. Using the group structure \ref{uiwrewd} of $L \cA\rightarrow \cA$ and the fact that the groupoid
$\cA$ is {\'e}tale it follows that $1\equiv t^n:T\rightarrow L\cA$ after replacing $T$ by some neighbourhood of $*$, if necessary  (see the proof of Lemma \ref{ewhjd} for a similar argument). It follows that $t^n:T\rightarrow L\cA$ has a natural lift $\tilde t^n$ given by an $H$-translate of the identity map such that $\sigma^n=\tilde t^n(*)$.

It remains to find the $n$`th root $\tilde t$ of $\tilde t^n$.
We now consider the diagram
$$\xymatrix{\ker(\dots)^n\ar[d]\ar[r]&H\ar[d]\ar[r]^{(\dots)^n}\ar[r]&H\ar[d]\\\ker(\dots)^n\ar[r]&\tilde \cG\times_{L\cA}T\ar[r]^{(\dots)^n}_{c}&\tilde \cG\times_{L\cA}T}\ .$$
The map
$c:\tilde \cG\times_{L\cA}T\rightarrow \tilde \cG\times_{L\cA}T$ is \'etale.
Therefore, after replacing $T$ by some neighbourhood of $*$ again, the datum of $\sigma$ and $\tilde t^n$ 
give the unique lift $\tilde t$.\hB  

\subsubsection{}

For smooth gerbes with band $U(1)$ on orbifolds the analog of Lemma \ref{flat} was shown e.g. in \cite{math.KT/0505267} or \cite{math.AT/0512658}. The argument in these papers uses the existence of a geometric structure (connection and curving) on the gerbe $G$. This geometry naturally induces a connection on the $U(1)$-principal bundle $\tilde G\rightarrow LX$. By a calculation the curvature of this connection vanishes. This gives the reduction of structure groups.

\subsubsection{}

Let $g:Y\rightarrow X$ be a map of topological stacks and $f:G\rightarrow X$ be a topological gerbe with band $H$ over $X$. We consider a $2$-cartesian diagram
$$\xymatrix{K\ar[d]\ar[r]&G\ar[d]\\Y\ar[r]^g\ar@{:>}[ur]&X}\ .$$
\begin{lem}\label{functorall}
We have a $2$-cartesian diagram
\begin{equation}\label{wueidiuwed}
\xymatrix{\widetilde{ K}\ar[d]\ar[r]^\phi&\tilde G\ar[d]\\
LY\ar[r]^{Lg}\ar@{:>}[ur]&LX}\ .
\end{equation}
Under the assumptions of Lemma \ref{flat}
this diagrams refines to
a $2$-cartesian diagram
$$\xymatrix{\tilde K^\delta\ar[d]\ar[r]&\tilde G^\delta\ar[d]\\
LY\ar[r]^{Lg}\ar@{:>}[ur]&LX}\ .$$
\end{lem}
\proof
We get the square (\ref{wueidiuwed}) from the functoriality part of  Proposition
\ref{ewiudiwed}. Since the vertical maps are $H$-principal bundles it is automatically $2$-cartesian.
The second statement easily follows from Lemma \ref{flat}. \hB

% {\tiny We argue on the level of topological groupoids.
% We have the following diagram
% $$\xymatrix{(LK)^0\ar[r]\ar[ddr]\ar[d]&K^1\ar[d]\ar[r]&G^1\ar[d]\\(LY)^0\ar[ddr]\ar[r]&Y^1\ar[r]&X^1\\&(LG)^0\ar[d]\ar[uur]&\\&(LX)^0\ar[uur]&}\ .$$
% In view of (\ref{deftg}) we must show that the left lower square is $2$-cartesian.
% But thus is true since the remaining three squares are $2$-cartesian.
% 
% Using the fact that $\phi$ is compatible with the group structures 
% it is easy to check the second assertion.}\footnote{\bf Dieses Argument geht so nicht. Am besten ist es, die Eigenschaft, daÃ der Loop funktor cartesische Diagramme erhÃ¤lt, auszunutzen. Wenn man $\tilde G$ durch $LG$ ersetzen wÃ¼rde, dann wÃ¤re vieles einfacher.} 

\hB

\subsection{The holonomy of $\tilde G^\delta$}\label{cycl1}

\subsubsection{}
Let $G\to X$ be a topological gerbe with band $U(1)$ over an orbispace $X$. In \ref{twto} we constructed a $U(1)^\delta$-principal bundle $G^\delta\to LX$.
It is an instructive exercise  to calculate the holonomy of this bundle in terms of the Dixmier-Douady invariant $d\in H^3(X;\Z)$ of the gerbe $G\to X$.
In the following we consider a special but typical case of this problem.

\subsubsection{}

We consider a $U(1)$-principal bundle
$\pi:E\to B$ in orbispaces and a topological gerbe $f:G\to E$ with band $U(1)$.
Let $h:LB\to U(1)^\delta$ be the function associated to the principal bundle $E\to B$
as in Lemma \ref{ewhjd} and define $LB_1:=h^{-1}(1)$. Then by Lemma \ref{hnsdkjc} we have an induced $U(1)$-principal bundle $L\pi:LE\to LB_1$. The holonomy of the bundle
$\tilde G^\delta\to LE$ along the fibres of $L\pi$ gives rise to a function
$$g:LB_1\to U(1)^\delta$$ (see \ref{hjadsca} for a precise construction).

The gerbe $f:G\to B$ is classified by a Dixmier-Douady class $d\in H^3(E;\Z)$.
Let $\pi_!:H^3(E;\Z)\to H^2(B;\Z)$ be the integration map. According to 
\ref{bxncmy} the class $\pi_!(d)\in H^2(B;\Z)$  gives rise to a function
$$\overline{\pi_!(d)}:LB\to U(1)^\delta\ .$$
The main result of the present subsection is the following proposition.
\begin{prop}\label{nhedvhjeada}
We have the equality
$$g=\overline{\pi_!(d)}_{|LB_1}\ .$$
\end{prop}

\subsubsection{}\label{hjadsca}

Here is the precise construction of the function $g:LB_1\to U(1)^\delta$.
Let $T$ be a space and $T\to LB_{1}$ be a map. The pull-back
$$\xymatrix{W\ar[r]\ar[d]&\tilde G^\delta\ar[d]\\S\ar@{:>}[ur]\ar[d]\ar[r]&LE\ar[d]\\ T\ar@{:>}[ur]\ar[r]&LB_{1}}$$ defines a  $U(1)$-principal bundle $S\to T$ and a $U(1)^\delta$-principal bundle $W\to S$. We chose an open  covering $(T_\alpha\to T)_{\alpha\in I}$ such that for all $\alpha\in I$ there exists a section
$$\xymatrix{&S\ar[d]\\T_\alpha\ar[r]\ar@{.>}[ur]^{s_\alpha}&T}\ .$$
The section $s_\alpha$  gives rise to a map
$T_\alpha\times \R\to S$ by $(t,x)\mapsto s_{\alpha}(t)x$, where $\R$ acts on $S$ via the covering $\R\to U(1)$. We can now (after refining the covering $(T_\alpha\to T)$ if necessary) choose a lift
$$\xymatrix{&W\ar[d]\\T_\alpha\times \R\ar@{.>}[ur]^{w_\alpha}\ar[r]&S}\ .$$
Then we define a map $g_{T_\alpha}:T_\alpha\to U(1)^\delta$ such that
$w_\alpha(t,0)=w_\alpha(t,1)g_{T_\alpha}(t)$.
Observe that $g_{T_\alpha}$ does not depend on the choices of $s_\alpha$ and $w_\alpha$.
One easily checks that the family  of maps $(g_{T_\alpha})_{\alpha\in I}$ determines a map
$g_T:T\to U(1)^\delta$ which depends functorially on $T\to LB_{1}$. It therefore defines a map
$g:LB_{1}\to U(1)^\delta$.

\subsubsection{}\label{qukendbma}
We now turn to the actual proof of  Proposition \ref{nhedvhjeada}.
We first consider a special case.
Let $\Gamma$ be a finite cyclic group which we write additively. We let $\Gamma$ act trivially on $U(1)$ and consider the orbispace $E:=[U(1)/\Gamma]$. The projection $U(1)\to *$ induces a $U(1)$-principal bundle $\pi:E\to B:=[*/\Gamma]$. We calculate $H^3(E;\Z)$ using the Kuenneth formula
and the product decomposition $E=U(1)\times B$. Note that $H^*(B;\Z)\cong H^*(\Gamma;\Z)$. In particular we have $H^3(B;\Z)\cong 0$ and a canonical isomorphism $H^2(B;\Z)\cong \hat \Gamma$ (see \ref{ejfwws}). It follows that 
$$H^3(E;\Z)\cong H^1(U(1);\Z)\otimes H^2(B;\Z)\cong \hat \Gamma$$
using the canonical orientation 
$H^1(U(1);\Z)\cong \Z$  of $U(1)$.

\subsubsection{}
 The group $H^3(E;\Z)$ classifies topological $U(1)$-gerbes over $E$.
In the following we present a construction which associates to every character
$\phi\in \hat \Gamma$ a $U(1)$-gerbe $G_\phi\to E$. We construct these gerbes in terms of representing groupoids.

The canonical covering $\R\to U(1)$ induces an atlas
$\R \rightarrow E$. The corresponding topological groupoid is the action groupoid for the action of $\Z\times \Gamma$ on $\R$ by $(n,\gamma)t:=t+n$. It is given  by
\begin{equation}\label{jdeikkike}
\R\times \Z\times \Gamma \Rightarrow \R
\end{equation}
with range $r(t,n,\gamma):=t+n$, source $s(t,n,\gamma):=t$, and
the composition $(t+m,n,\gamma)\bullet(t,m,\gamma^\prime):=(t,n+m,\gamma+\gamma^\prime)$.

The character $\phi\in \hat \Gamma$ determines a $U(1)$-central extension
\begin{equation}\label{ldedjedkeded}
0\to U(1)\to \widehat{\Z\times \Gamma}_\phi\to \Z\times \Gamma\to 0\ .
\end{equation}
If we identify $\widehat{\Z\times \Gamma}_\phi\cong \Z\times  \Gamma\times U(1)$ as sets, then the multiplication is given by
$(n,\gamma,z)(n^\prime,\gamma^\prime,z^\prime)=(n+n^\prime,\gamma+\gamma^\prime,\phi(\gamma)^{n^\prime} z z^\prime)$. This central extension acts on $\R$ via its projection
$\widehat{\Z\times \Gamma}_\phi\to \Z\times \Gamma$, $(n,\gamma,z)\mapsto (n,\gamma)$. The gerbe
$G_\phi\to E$ is then given by
$$[\R/\widehat{\Z\times \Gamma}_\phi]\to [\R/\Z\times \Gamma]\ .$$ In terms of groupoids, $G_\phi$ is 
given as the  $U(1)$-central extension of the groupoid (\ref{jdeikkike}) which on the level of morphisms is the trivial $U(1)$-bundle
$$\R\times \Z\times \Gamma\times U(1)\rightarrow \R\times \Z\times \Gamma\ ,$$ whose source and range maps are
$$s(t,n,\gamma,z):=t\ , \quad r(t,n,\gamma,z):=t+n\ ,$$
and whose composition is given by 
$$(t+m,n,\gamma,z^\prime)(t,m,\gamma^\prime,z):=(t,n+m,\gamma+\gamma^\prime,\phi(\gamma)^m z^\prime z)\ .$$

\subsubsection{}

We now calculate the bundle $\tilde G_{\phi}^\delta\rightarrow LE$. First of all note that $$LE \cong [\Gamma\times U(1)/\Gamma]\ ,$$ where $\Gamma$ acts trivially on $\Gamma\times U(1)$. The map $\Gamma\times \R\to \Gamma\times U(1)$ gives an atlas of $LE$.
The associated groupoid is the action groupoid of the action of $\Z\times \Gamma$ on
$\Gamma\times \R$ by $(n,\gamma)(\sigma,t)=(\sigma,t+n)$.
It is given by  
$$\Gamma\times \R\times \Z \times \Gamma\Rightarrow \Gamma\times \R$$
with range and source given by
$$r(\sigma,t,n,\gamma):=(\sigma,t+n)\ , \quad s(\sigma,t,n,\gamma):=(\sigma,t)\ ,$$
and with the composition
$$(\sigma,t+m,n,\gamma)\circ (\sigma,t,m,\gamma^\prime):=(\sigma,t,n+m,\gamma+\gamma^\prime)\ .$$
We can now read off a groupoid
 presentation of the   $U(1)^\delta$-principal bundle
$\tilde G_{\phi}^\delta\rightarrow LE$. It is
is presented by the $U(1)^\delta$-principal bundle in groupoids
$$\xymatrix{ 
\Gamma\times \R\times U(1)^\delta  \times \Z\times \Gamma\ar[d]\ar@{=>}[r]&\Gamma\times \R\times U(1)^\delta\ar[d]\\
\Gamma\times \R\times \Z\times \Gamma\ar@{=>}[r]&\Gamma\times \R }\ .$$
The range and source maps in the upper horizontal line are given by 
$$r(\sigma,t,z,n,\gamma):=(\sigma,t+n,\phi(\sigma)^nz)\ , \quad s(\sigma,t,z,n,\gamma):=(\sigma,t,z)\ ,$$
and with the composition
$$(\sigma,t+m,\phi(\sigma)^mz,n,\gamma)\circ (\sigma,t,z,n,\gamma^\prime):=(\sigma,t,z,n+m,\gamma+\gamma^\prime)\ .$$
In particular, the holonomy of $\tilde G_{\phi}^\delta$ along the fibre of $LE$ over 
$[\{\sigma\}/\Gamma]$ is given by $\phi(\sigma)$.

\subsubsection{}

In our example we  have $LB_1=[\Gamma/\Gamma]=LB$, where $\Gamma$ acts trivially on itself.
The function $g_\phi:LB_1\to U(1)^\delta$, which measures the holonomy of $\tilde G_\phi\to LE$ along the fibres of $LE\to LB_1$, is given by the calculation above by 
\begin{equation}\label{qewgd}g_\phi=\phi:\Gamma\to U(1)^\delta\ .\end{equation} 
By the discussion \ref{qukendbma} the character
$\phi$ gives rise to a class $d_\phi\in H^3(E;\Z)$ such that
$$\pi_!(d_\phi)=\phi$$ (using the isomorphism $\hat \Gamma\cong  H^2(B;\Z)$).
Furthermore we have
$$\overline{\pi_!(d_\phi)}=\phi:\Gamma\to U(1)^\delta\ .$$

\subsubsection{}

In order to finish the proof of Proposition \ref{nhedvhjeada} in the special case  we must show that $d_\phi$ is the Dixmier-Douady class $d(G_\phi)$ of $G_\phi$.
We will use the following two general facts:
\begin{enumerate}
\item
Let $1\to U(1)\to \hat G\to G\to 1$ be a $U(1)$-central extension of a discrete group $G$ classified by $e\in \Ext(G;U(1)):=H^2(G;U(1))$. Furthermore, let
$\delta:H^2(G;U(1))\to H^3(G;\Z)$ be the boundary operator in the Bockstein sequence 
in group cohomology associated to the exact sequence of coefficients
$0\to \Z\to \R\to U(1)\to 0$.
Then the Dixmier-Doudady class of the gerbe
$[*/\hat G]\to [*/G]$ is given by the image of 
$\delta(e)\in H^3(G;\Z)$ under the isomorphism $H^3(G;\Z)\cong H^3([*/G];\Z)$.
\item
Let $\phi:G\to U(1)$ be a character of a finite group $G$.
It gives rise to a class $\phi\in H^1(G;U(1))$ and an extension
 $1\to U(1)\to \widehat{\Z\times G}\to \Z\times G\to 1$. We can identify $\widehat{\Z\times G}\cong \Z\times G\times U(1)$ as sets. Its multipication  is then  given by $(n,g,z)(n^\prime,g^\prime,z^\prime)=(n+n^\prime,gg^\prime,\phi(g)^{n^\prime}zz^\prime)$.
The class $e\in \Ext(\Z\times G;U(1))$ of the extension
is then given by image of 
$\id_\Z\times \phi\in H^1(\Z;\Z)\times H^1(G;U(1))$
under the product
$\times:H^1(\Z;\Z)\times H^1(G;U(1))\to H^2(\Z\times G;U(1))$,
where $\id_\Z\in H^1(\Z;\Z)$ is the identity homomorphism.
\end{enumerate}

We now specialize these facts to the present situation.
The Kuenneth formula gives an isomorphism
\begin{equation}\label{dkjedwwww}
\Ext(\Z\times \Gamma;U(1)):=H^2(\Z\times \Gamma;U(1))\cong H^1(\Z;\Z)\otimes H^1(\Gamma;U(1))\cong 
\Z\otimes \hat \Gamma\cong \hat \Gamma\ ,\end{equation}
where we use the generator $\id_\Z\in H^1(\Z;\Z)$ in order to identify   
$H^1(\Z;\Z)\cong \Z$.
The class $e_\phi\in \Ext(\Z\times\Gamma;U(1))$ of the extension  (\ref{ldedjedkeded})  corresponds under this isomorphism to $\phi\in \hat \Gamma$ (by (2)).

% Let $$\delta:H^2(\Z\times\Gamma;U(1))\to H^3(\Z\times \Gamma;\Z)$$
% denote the boundary map of the Bockstein sequence associated to the sequence of coefficients  $$0\to \Z\to \R\to U(1)\to 0\ .$$
By (1) the Dixmier-Douady class  $d(\phi)\in  H^3([*/\Z\times \Gamma];\Z)$ of the gerbe
$$[*/\widehat{\Z\times \Gamma}_\phi]\to [*/\Z\times \Gamma]$$ corresponds to $$\delta(e_\phi)\in H^3(\Z\times \Gamma;\Z)$$    under the identification
$$H^3([*/\Z\times \Gamma];\Z)\cong H^3(\Z\times \Gamma;\Z)\ .$$

\subsubsection{}

Let $p: B=[\R/\Z\times \Gamma]\to [*/\Z\times \Gamma]$ be the canonical projection.
Then we have 
$d(G_\phi)=p^*d(\phi)$.
We now observe that the following diagram commutes,
$$\xymatrix{d(G_\phi)\ar@/_2cm/[ddd]&d(\phi)\ar@/_-0.4cm/[l]&&e_\phi\ar@/_-0.5cm/[ll]\ar@/_-2.5cm/[dddl]^{(\ref{dkjedwwww})}\\H^3([\R/\Z\times \Gamma];\Z)\ar[d]^{\pi_!}&H^3([*/\Z\times \R];\Z)\ar[l]^{p^*}_\cong &H^3(\Z\times \Gamma;\Z)\ar[l]^\cong&H^2(\Z\times\Gamma;U(1))\ar[d]^\cong\ar[l]_\cong^\delta\\H^2([*/\Gamma];\Z)&H^2(\Gamma;\Z)\ar[l]^\cong&H^1(\Gamma;U(1))\ar[l]^\delta&H^1(\Z;\Z)\otimes H^1(\Gamma;U(1))\ar[l]^\cong\\\pi_!(d_\phi)&&\phi\ar@/^-0.5cm/[ll]^{definition\: of\: d_\phi}& }\ ,$$
and that the elements are mapped as indicated.

% \subsubsection{}
% The set of isomorphism classes of $U(1)$-gerbes over an orbispace  $E$ forms a group with respect to the "tensor" product of gerbes, and the Dixmier-Douady class identifies this group with the 
% additive group $H^3(E;\Z)$. 
% By construction for $k\in \Z$  we have $G_{\phi^k}\cong G_\phi^{\otimes k}$. Therefore  $$d(G_{\phi^k}) = k d(G_\phi)\ .$$
% 
% It follows from (\ref{qewgd}) the $\phi\not=1$ implies that $G_\phi$ is not trivial. In particular, if $\phi$ generates the group $\hat \Gamma$, then $d(G_{\phi})$ generates
% the group $H^3(E;\Z)$. In this case $d_\phi\in H^3(E;\Z)$ is another generator.
%  

\subsubsection{}

We show how the general case of Proposition \ref{nhedvhjeada} can be reduced to the special case discussed above.
The constructions of $g$ and $\overline{\pi_!(d)}_{|LB_1}$ are natural with respect to pull-back. Therefore in order to verify Proposition \ref{nhedvhjeada}  it suffices to show the desired equality over each point in $LB$ seperately. As in \ref{bxncmy} a point $u\in LB$ is given by a 
point $p\in B$ and an element $\gamma\in \Aut(p)$ (in the present subsection we omit the subscript $u$ in order to simplify the notation).  Let $\Gamma\subset  \Aut(p)$ be the cyclic group generated by $\gamma$ and $\chi\in \hat \Gamma$ be the character by which $\Gamma$  acts on the fibre $\pi^{-1}(p)$. Note that
\begin{equation}\label{wzehbd}
\chi(\gamma)=h(u)\ .
\end{equation}
We get a cartesian diagram
\begin{equation}\label{weikdnasm}\xymatrix{v^*G\ar[d]\ar[r]&G\ar[d]^f\\ [U(1)/_\chi\Gamma]\ar[r]^v\ar[d]^q&E\ar[d]^\pi\\[*/\Gamma]\ar[r]^{\tilde u}&B}
\ .\end{equation}
such that $L\tilde u(\gamma)=u$, 
where we consider $\gamma\in [\Gamma/\Gamma]\cong L[*/\Gamma]$. 
In particular, $v^*d$ is a Dixmier-Douady class of the gerbe $v^*G\to [U(1)/_\chi\Gamma]$ and  we have 
$$\overline{\pi_!(d)}(u)=\overline{q_!v^*(d)}(\gamma)\ .$$ 
Observe that $L[*/\Gamma]_1=[\ker(\chi)/\Gamma]$.
Let $g_{v^*G}:L[*/\Gamma]_1\to U(1)^\delta$ denote the function (\ref{hjadsca}) which measures the holonomy of $\widetilde{v^*G}^\delta\to L[U(1)/_\chi\Gamma]$ along the fibres of $q$. If $u\in LB_1$, then by \ref{wzehbd} we have $\chi(\gamma)=1$ and 
$$g(u)=g_{v^*G}(\gamma)\ .$$
The equation $$\overline{\pi_!(d)}(u)=g(u)$$ now follows from the equation
$$ g_{v^*G}(\gamma)=\overline{q_!v^*(d)}(\gamma)$$
which was already shown above. \hB

\section{Delocalized Cohomology of orbispaces and orbifolds}\label{zuwlklklklkdee}

\subsection{Definition of delocalized twisted cohomology}\label{duedeeeww}

\subsubsection{}

A topological stack $X$ gives rise to a site $\Site(X)=\bX$.
The underlying category of $\bX$ is the subcategory of $\Top/X$ of maps
$(U\to X)$ which are representable and have local sections.
The covering families $(U_i\to U)$ are families of maps $U_i\to U$ in $\bX$ which have local sections\footnote{A map of topological spaces $f:V\to W$ has local sections if for every $w\in f(V)$ there exists an open  neighbourhood $W_w\subseteq W$ and a map $\sigma:W_w\to V$ such that $\id_{W_w}=f\circ \sigma$} and are such that $\sqcup_{i} U_i\to U$ is surjective. 
One can actually restrict to covering families by open subsets without changing the induced topology (the argument is similar as for  \cite[Lemma 2.47]{bss}).
If $X$ is a space, then the small site $(X)$ of $X$ is the category of open subsets of $X$
with the usual notion of covering families.

\subsubsection{}
To the site $\bX$ we associate categories of presheaves
and sheaves $\Pr\bX$ and $\Sh\bX$ in the usual way. A map $p:X\to Y$ of topological stacks induces
a pair of adjoint functors
$$p^*:\Sh\bY\Leftrightarrow \Sh\bX:p_*\ .$$
We use this framework of sheaf theory on topological stacks in order to define the delocalized cohomology of an orbispace twisted by a gerbe.

For details of the sheaf theory we refer to \cite{bss} and \cite{bssf}.

For a site $\bX$ let $i:\Sh \bX\to \Pr\bX$ denote the canonical embedding
of the category of presheaves into the category of sheaves, and let $i^\sharp:\Pr\bX\to \Sh\bX$ denote its left-adjoint, the sheafification functor. We use the same symbols in order to denote the restriction of these functors to the categories $\Pr_\Ab\bX$ and $\Sh_\Ab\bX$ of presheaves and sheaves of abelian groups.

\subsubsection{}\label{diwued}
Let $H$ be a topological abelian group. We assume that $\Tors_n(H)\subseteq H$ is discrete for all $n\in \nat$ (see \ref{ejkjkjkked666}). Let $H^\delta$ denote the group $H$ with the discrete topology.
Furthermore, let $Z$ be a discrete  abelian group and $\lambda:H^\delta\rightarrow \Aut(Z)$ be a homomorphism. 

\subsubsection{}\label{grs}

Let $P\rightarrow X$ be the underlying map of stacks of an $H^\delta$-principal bundle  over a topological stack $X$.  If $(U\to X)\in \bX$, then $U\times_XP\rightarrow U$ is an ordinary  $H^\delta$-principal bundle.
We define the abelian group $\cZ_{P,\lambda}(U)$ to be the group of  continuous  sections of the associated bundle
$(U\times_XP)\times_{H^\delta,\lambda}Z\rightarrow U$ under pointwise multiplication.
If $(U^\prime\to X)\rightarrow (U\to X)$ is a morphism in $\bX$, then we have an induced morphism $U^\prime \times_X P\to U\times_XP$ of $H^\delta$-principal bundles over $U^\prime\to U$. It induces a homomorphism $\cZ_{P,\lambda}(U)\to \cZ_{P,\lambda}(U^\prime)$. 
In this way obtain a presheaf of abelian groups $\cZ_{P,\lambda}\in \Pr_\Ab\bX$,
$U\mapsto \cZ_{P,\lambda}(U)$. 
Note that $\cZ_{P,\lambda}$ is actually a sheaf, i.e. we have $\cZ_{P,\lambda}\in \Sh_\Ab\bX$.

\subsubsection{}\label{gggdqiowddd}

Let $f:G\rightarrow X$ be a gerbe with band $H$ over an orbispace $X$. Then by Lemma \ref{flat} we have the $H^\delta$-principal bundle $\tilde G^\delta\rightarrow LX$. By \ref{grs} it gives rise to a the presheaf
$\cZ_{\tilde G^\delta,\lambda}\in \Pr_\Ab\bLX$.  

\subsubsection{}
We define a gerbe $f_L:G_L\rightarrow LX$ with band $H$ as the pull-back of the gerbe $f:G\to X$ along the canonical map $i:LX\to X$ (see \ref{kk77qghhgsghvqhgsq}). We have a diagram
$$\xymatrix{{*}&G_L\ar[r]\ar[l]^p\ar[d]^{f_L}&G\ar[d]^f\\&LX\ar@{:>}[ur]\ar[r]^i&X}\ .$$
We  consider $f_L^*\cZ_{\tilde G^\delta,\lambda}\in \Sh_\Ab\bGL$.

\subsubsection{}

Let $\ev:=\ev_{*\to *}:\Sh_\Ab\Site(*)\to\Ab$ be the functor, which evaluates a sheaf of abelian groups on $\Site(*)$ at the object
$(*\to *)\in \Site(*)$.
\begin{lem}
The functor $\ev: \Sh_\Ab\Site(*)\to\Ab$ is exact.
\end{lem}
\proof
A basic observation lying at the heart of sheaf theory is that evaluation functors  are not exact in general. Therefore, a proof of exactness of the evaluation $\ev$ is required.
First note that $\Site(*)$ is the big site of $*$ which can be identified with the category of all topological spaces.
Every non-empty collection of non-empty spaces is a covering family of $*$.

The small site $(*)$ of $*$ has one object $*\to *$.  In \cite{bssf} (see also \cite[Prop. 2.46]{bss}, the arguments works equally well in the smooth and topological contexts)  we have seen that the restriction functor $\nu_*:\Sh\Site(*)\to \Sh(*)$ is exact. 
Let $\tilde\ev:\Sh_\Ab(*)\to \Ab$ denote the corresponding evaluation functor. It is actually an isomorphism of categories, and in particular exact.
We have $\tilde \ev\circ \nu_*\cong \ev$. We see that  $\ev$ is exact, since it is a composition of exact functors.
\hB

\subsubsection{}

The functor $p_*:\Sh_\Ab(\bGL)\to \Sh_\Ab\Site(*)$ is left-exact and thus admits right-derived functors
$$Rp_*:D^+( \Sh_\Ab\bGL)\to  D^+(\Sh_\Ab\Site(*))$$
between the lower bounded derived categories. The functor
$ev:\Sh_\Ab\Site(*)\to \Ab$ is exact and thus descends to the lower-bounded derived catgeories.

\begin{ddd}\label{defdelr}
We define the delocalized $G$-twisted cohomology of $X$ with coefficients in $(Z,\lambda)$ by
$$H^*_{deloc}(X;G,(Z,\lambda)):=H^*(\ev\circ Rp_*  (f^*_L\cZ_{\tilde G^\delta,\lambda}
))\ .$$
\end{ddd}

\subsubsection{}\label{ldefg}

The most important example for us is the case where $Z:=\C^\delta$ and $H:=U(1)$ with $\lambda:H^\delta\rightarrow Z\rightarrow \End(Z)$ being the obvious embedding $U(1)^\delta\to \End(\C^\delta)$. In this case we will denote the sheaf $\C_{\tilde G^\delta,\lambda}$ by $\cL$ or $\cL_G$, if a reference to $G$ is necessary.
\begin{ddd}\label{fdifdfdmhjjhhjjh}
The $G$-twisted complex delocalized cohomology of $X$ is defined by 
$$H^*_{deloc}(X;G):=H^*_{deloc}(X;G,\cL)\ .$$
\end{ddd}

\subsubsection{}

Another example related to $Spin$-structures is the case where
$Z:=\Z$, $H:=\Z^*=\{1,-1\}$, and $\lambda:\Z^*\to \End(\Z)$ is again the canonical embedding.

\subsubsection{}\label{qiodiowqodqwdqw}

We now discuss the functial behaviour of the delocalized twisted cohomology.
We defined the sheaf $\cZ_{\tilde G^\delta,\lambda}$ down on $LX$ in order to connect with usual conventions in the literature on inner local systems and twisted torsion, and in order to have the formula (\ref{hgghsd4557887812}) below.
This construction depends on descending the $H^\delta$-bundle $LG\to G_L$
to the bundle $\tilde G\to LX$. The quite  complicated construction was carries out in 
Proposition \ref{ewiudiwed}. In the definition of twisted cohomology we then use the pull-back $f_L^*\cZ_{\tilde G^\delta,\lambda}$.  

It would be much more natural to construct the sheaf
$\tilde \cZ_{LG^\delta,\lambda}:=f_L^*\cZ_{\tilde G^\delta,\lambda}$ directly starting from the $H^\delta$-principal bundle $LG^\delta\to G_L$.
We can proceed as in the definition of $\cZ_{\tilde G^\delta,\lambda}$. For an object $(U\to G_L)\in \bGL$ we define $\tilde \cZ_{LG^\delta;\lambda}(U)\in \Ab$ as the group of continuous sections of $(U\times_{G_L}LG^\delta)\times_{H^\delta,\lambda} Z$ under pointwise multiplication. For a morphism
$U^\prime\to U$ we then have a natural homomorphism $\tilde \cZ_{LG^\delta;\lambda}(U)\to \tilde \cZ_{LG^\delta;\lambda}(U^\prime)$ induced by a corresponding morphism of principal bundles over $U^\prime\to U$.

We have a canonical isomorphism
$$H_{deloc}(X;G,(Z,\lambda))\cong H^*(\ev\circ Rp_*(\tilde \cZ_{LG^\delta,\lambda}))\ .$$
In the case $H=U(1)$ and $Z=\C^\delta$ we set $\tilde \cZ_{LG^\delta,\lambda}:=\tilde \cL$.
\subsubsection{}

We consider a $2$-cartesian diagram
$$\xymatrix{G^\prime\ar[d]^{f^\prime}\ar[r]^h&G\ar[d]^f\\X^\prime\ar@{:>}[ur]\ar[r]^g&X}\ ,$$ where $g$ is a map of orbispaces, i.e. a representable map.

\begin{lem}\label{zuzuzwdwdwdw}
We have a canonical functorial map
$$(g,h)^*:H^*_{deloc}(X;G)\to H^*_{deloc}(X^\prime;G^\prime)\ .$$
\end{lem} 
\proof
Since the loop functor preserves two-cartesian diagrams
we get an induced $2$-cartesian diagram
\begin{equation}\label{weuuuwk88z7fwef}
\xymatrix{LG^{\prime\delta}\ar[d]\ar[r]^{Lh}&LG\ar[d]\\G^\prime_L\ar@{:>}[ur]\ar[d]^{f^\prime_L}\ar[r]^{h_L}&G_L\ar[d]^{f_L}\\LX^\prime\ar@{:>}[ur]\ar[r]^{Lg}&LX}\ .
\end{equation}
Let $\tilde \cL=\tilde \cZ_{LG^\delta,\lambda}\in \Sh_\Ab \bGL$ and $\tilde \cL^\prime:=\cZ_{LG^{\prime\delta},\lambda}\Sh_\Ab\bGL^\prime$ denote the sheaves of abelian groups  associated to
$G$ and $G^\prime$ and $(Z,\lambda)$ as in \ref{qiodiowqodqwdqw}. The diagram 
(\ref{weuuuwk88z7fwef}) induces an isomorphism
\begin{equation}\label{shjssa}
h_L^*\tilde \cL\stackrel{\sim}{\to} \tilde \cL^\prime
 \end{equation} of sheaves on $G^\prime_L$.
We now consider the diagram
$$\xymatrix{G^\prime_L\ar[dr]^{p^\prime}\ar[rr]^{h_L}&&G_L\ar[dl]_{p}\\&{*}}\ .$$  The unit  
$\id\to R(h_L)_*\circ h_L^*$  of the adjoint pair 
$$h_L^*:D^+(\Sh_\Ab\bGL)\Leftrightarrow D^+(\Sh_\Ab \bGL^\prime):R(h_L)_*$$  induces a natural transformation
\begin{equation}\label{kjnjwdwwwww}
Rp_*\to Rp_*\circ R(h_L)_*\circ  h_L^*:D^+(\Sh_\Ab \bGL)\to D^+(\Sh_\Ab\Site(*))\ .\end{equation}
Since $p\circ h_L=p^\prime$ and $Rp_*\circ R(h_L)_*\cong R(p\circ h_L)_*$ (see \cite{bssf} and also \cite[Lemma 2.26]{bss} for an argument in the smooth case)  we have 
an isomorphism 
$$Rp_*\circ R(h_L)_*\cong R(p\circ h_L)_*\cong Rp^\prime_*\ .$$
We insert this into (\ref{kjnjwdwwwww})
and get the natural transformation
\begin{equation}\label{wjhzkjkcowww}
Rp_*\to Rp^\prime_*\circ  h_L^*:D^+(\Sh_\Ab \bGL)\to D^+(\Sh_\Ab\Site(*))
\ .\end{equation}
We define
$$(g,h)^*:H^*\circ \ev\circ Rp_* (\tilde \cL)\stackrel{(\ref{wjhzkjkcowww})}{\to} H^*\circ \ev\circ Rp^\prime_*\circ h_L^*  (\tilde\cL)\stackrel{(\ref{shjssa})}{\cong} H^*\circ \ev\circ Rp^\prime_*(\tilde \cL^\prime)\ .$$
We leave it to the reader to write out the argument for functoriality. The basic input is the functoriality of the units for a composition $f\circ g$ which can be expressed as the commutativity of  
$$\xymatrix{\id\ar[r]\ar@/_-0.5cm/[rrr]&Rf_*\circ f^*\ar[r]&Rf_*\circ Rg_*\circ g^*\circ f^*\ar[r]^{\cong}&R(f\circ g)_*\circ (f\circ g)^*}$$
(see \cite{bssf} for a proof).
\hB

\subsubsection{}\label{umr}

>From now on we consider the case $H:=U(1)$ and $Z:=\C^\delta$.
We can decompose
$p=q\circ f_L$, where $q:LX\rightarrow *$.
Since $f_L$ has local sections we have an isomorphism $$Rp_*\cong Rq_*\circ R(f_L)_*$$ by \cite[2.26]{bss}.
We have
$$Rp_*\circ f_L^*(\cL)\cong Rq_*\circ R(f_L)_* \circ f_L^*(\cL)$$
and the projection formula (see \cite{bssf})
\begin{equation}\label{hgghsd4557887812}
R(f_L)_*(\tilde \cL)\cong R(f_L)_* \circ f_L^*(\cL)\cong R(f_L)_* (\C_{\bG_L})\otimes_\C \cL\ .
\end{equation}
 Therefore we can write
$$H^*_{deloc}(X;G)\cong H^*(\ev\circ Rq_*\circ (R(f_L)_* (\C_{\bG_I})\otimes_\C \cL))\ .$$

\subsection{Twisted de Rham cohomology}

\subsubsection{}

The theory developed in the Sections \ref{ttzwtzetzwiu}, \ref{wsdwkjdwds}, \ref{twto} and \ref{duedeeeww} has a counterpart in the world of stacks in smooth manifolds though there is one essential difference. The map $LX\to X$ is not representable as a map of stacks in smooth manifolds.
Therefore the proof of the fact that $LX$ is a smooth stack is quite different from the topological case\footnote{One could save our argument by introducing the notion of a smoothly representable map between stacks in smooth manifolds and showing that $LX\to X$ is smoothly representable. A map $X\to Y$ between stacks in smooth manifolds is called smootly representable, if the fibre product
$A\times_YX$ is a manifold for every {\em submersion} $A\to Y$.}. But note that we have not used representability of $LX\to X$ otherwise. 

 In  the following we explain the replacements which lead to a precisely analogous theory.
\begin{enumerate}
\item
The category of topological spaces $\Top$ is replaced by the category of $\Mf^\infty$ of smooth manifolds.
\item Stacks in topological spaces are replaced by stacks in manifolds.
\item The condition on a map of having local sections is replaced by the condition of being a submersion (following the conventions from algebraic geometry we will use the term "smooth" synonymously with "submersion").
\item Topological stacks are  replaced by smooth stacks. A stack in smooth manifolds $X$ is called smooth if it admits  an atlas $a:A\to X$, i.e.   a representable, surjective, and  submersive (which replaces the local section condition by the preceding point) map from a manifold $A$.
\item The notion of a topological groupoid is replaced by the notion of a Lie groupoid. In particular, we require that range and source maps are submersions.
 \item Orbispaces are replaced by orbifolds. A smooth stack is an orbifold if it admits an orbifold atlas. An orbifold atlas is an atlas which gives rise to a proper and \'etale groupoid in smooth manifolds. Since manifolds are locally compact and Hausdorff the conditions  "separated"  and "very proper"\footnote{The condition "very proper" is as in \ref{uiiuduiwed} with the difference that the cut-off function must be smooth.} hold automatically (see Lemma \ref{uduiqwdqw444}).
 \item The group $H$ in \ref{diwued} must be a Lie group. 
\item For a smooth stack $X$ the site $\bX$ is the subcategory of $\Mf^\infty/X$ of maps
$(U\to X)$ which are representable submersions. The covering families are families
$(U_i\to U)$ of submersions such that $\sqcup_iU_i\to U$ is surjective.
\end{enumerate}
One problem with the category $\Mf^\infty$ is that fibre products only exist under additional conditions (e.g. if one map is a submersion). We leave it to the interested reader to check that all fibre products used in Sections \ref{wsdwkjdwds}, \ref{twto} and \ref{duedeeeww}  exists in manifolds\footnote{with the exception that the construction of the simpler model of $L\cX\times_{\cX}\cG$ in \ref{hjqwhdjqwdq} needs different arguments since (\ref{hh32562e2398e83}) may not be  a transversal pull-back.}.

Let $X$ be an orbifold and $G\to X$ be a smooth gerbe with band $U(1)$.
Then by \ref{fdifdfdmhjjhhjjh} we have a well-defined twisted delocalized cohomology
$$H^*_{deloc}(X;G)\ .$$ 
The main goal of the present section is to calculate
this cohomology in terms of a twisted de Rham complex. 
This generalizes the main result of \cite{bss} from smooth manifold $X$ to orbifolds $X$.

\subsubsection{}

 The first goal of the present subsection is to define the de Rham complex associated to a locally constant sheaf of complex vector spaces on an orbifold in two equivalent (according to Lemma \ref{twopic}) ways. In the first picture we define a sheaf of de Rham complexes on the big site of the orbifold and then take its global sections. The second picture uses the calculus of differential forms on the orbifold itself. While the first picture belongs to the philosophy of the present paper this second definition is mainly used to compare with other constructions in the literature.

In the second part we apply this construction to the local system $\cL\in \Sh_\Ab\bLX$  associated to an $U(1)$-gerbe $G\to X$ on an orbifold.

\subsubsection{}

Consider a smooth stack $X$ in smooth manifolds. Let $\cE$ be a locally constant sheaf on $\bX$ of complex vector spaces.
If $(U\to X)\in \bX$, then $\cE_{|U}$ is the sheaf of parallel sections of a canonically determined complex vector bundle with flat connection $(E_U,\nabla^{E_U})$.
Let $\Omega^k(U,E_U)$ denote the space of global sections of
$\Lambda^k_\C T^*U\otimes E_U$. The de Rham differential $d_{dR}$ and the connection $\nabla^{E_U}$ together induce a differential
$d^{E_U}:\Omega^k(U,E_U)\rightarrow \Omega^{k+1}(U,E_U)$.
Observe that
$(\Omega^\cdot(U,{E_U}),d^{E_U})$ is a $(\Omega^\cdot(U),d_{dR})$-DG-module.

If $f:(U^\prime\to X)\rightarrow(U\to X)$ is a morphism in $\bX$, then we have a morphism of sheaves $f^*\cE_{|U}\rightarrow \cE_{|U^\prime}$. This induces a morphism of flat vector bundles
$f^*E_U\rightarrow E_{U^\prime}$ and finally a morphism of complexes
$(\Omega^\cdot(U,{E_U}),d^{E_U})\rightarrow (\Omega^\cdot(U^\prime,{E_{U^\prime}}),d^{E_{U^\prime}})$.

We define the sheaf $\Omega_X^\cdot(\cE)$ of $(\Omega^\cdot_X,d_{dR})$-DG-modules  which associates to $(U\to X)$ in $\bX$ the complex $(\Omega^\cdot(U,E_U),d^{E_U})$. 

\begin{lem}
$\cE\rightarrow \Omega_X^\cdot(\cE)$ is a flabby resolution
\end{lem}
\proof
This is shown by adapting the arguments of \cite[Sec. 3.1]{bss} to differential forms twisted by a flat vector bundle. \hB

\subsubsection{}\label{uieuhwuefuw}

Let $p:X\rightarrow *$ be the projection. If $F\in \Sh\bX$, then we define its global sections by $$\Gamma_X F:=\ev\circ p_*(F)\ .$$ 

\subsubsection{}\label{zueuedeke}
Assume now that $X$ is an orbifold.
The sheaf $\cE$ gives rise to a flat vector bundle $E\rightarrow X$ in the orbifold sense. We can consider the de Rham complex $\Omega(X,E)$ of $E$-valued forms on $X$ which are smooth in the orbifold sense. 

\begin{lem}\label{twopic}
We have a natural isomorphism $\Gamma_X\Omega_X^\cdot(\cE)\cong \Omega(X,E)$.
\end{lem}
\proof
We choose an orbifold atlas $A\rightarrow X$, i.e. $A$ is a smooth manifold, $A\rightarrow X$ is an atlas, and the smooth groupoid $A\times_XA\Rightarrow A$ is very proper, separated  and {\'e}tale.
By the definition of smooth forms in the orbifold sense we have the exact sequence
$$\xymatrix{0\ar[r]&\Omega(X,E)\ar[r]& \Omega(A,E_A)\ar[r]^{r^*-s^*\hspace{0.8cm}}&\Omega(A\times_XA,E_{A\times_XA})}\ .$$
The composition $A\to X\to *$ is clearly representable.
By \cite[Lemma 2.36]{bss} we have an exact sequence
$$\xymatrix{0\ar[r]&\Gamma_X \Omega^\cdot_X(\cE)\ar[r]& \Omega(A,E_A)\ar[r]^{r^*-s^*\hspace{0.8cm}}&\Omega(A\times_XA,E_{A\times_XA})}\ .$$
\hB

\subsubsection{}

In order to indicate that the local system $\cE$ is the initial datum and the vector bundle $E\to X$ is secondary we use the following notation.
\begin{ddd}\label{ewduiewdw33}
We define
$$\Omega^\cdot(X,\cE):=\Gamma_X\Omega^\cdot_X(\cE)\ .$$
It is a $\Omega^\cdot(X)$-DG-module. Its differential will be denoted by $d^\cE$.
\end{ddd}

\subsubsection{}
The twisted de Rham cohomology of $X$ with coefficients in $\cE$ depends on the choice of a closed form $\lambda\in \Omega^3(X)$. 
 Let $z$ be a formal variable of degree $2$. Then we 
form the complex $\Omega^\cdot(X,\cE)[[z]]_\lambda$ given by 
$$\Omega^\cdot(X,\cE)[[z]]\ ,\quad d_\lambda:=d^{\cE}+\lambda T\ ,$$
where $T:=\frac{d}{dz}$.

\begin{ddd}\label{dera}
The $\lambda$-twisted cohomology $H^*(X;\cE,\lambda)$ of $X$ with coefficients in $\cE$ is defined as the cohomology of the complex $\Omega^\cdot(X,\cE)[[z]]_\lambda$.
\end{ddd}

\subsubsection{}\label{dera1}

We can also define a sheaf
$\Omega^\cdot_X(\cE)[[z]]_\lambda$ of $(\cE,\lambda)$-twisted de Rham complexes on $\bX$  such that for $(U\stackrel{\phi}{\to} X)\in \bX$ we have
$\Omega^\cdot_X(\cE)[[z]]_\lambda(U):=\Omega^\cdot(U,E_U)[[z]]$ with the differential $d_{\phi^*\lambda}$. 
By Definition \ref{ewduiewdw33} have an isomorphism of complexes
$$\Omega^\cdot(X,\cE)[[z]]_\lambda \cong \Gamma_X\Omega^\cdot_X(\cE)[[z]]_\lambda \ .$$

\subsubsection{}\label{gensi}

We now take twists into account.
Let $X$ be an orbifold and  $f:G\rightarrow X$ be a smooth gerbe with band $U(1)$. Then we can form the orbifold of loops $LX\rightarrow X$ and the pull-back $f_L:G_L\rightarrow LX$ of the gerbe $f:G\to X$.
We choose an atlas $(A\rightarrow G_L)\in \bGL$. It gives rise to a simplicial object
$\bA^\cdot_{G_L}\in \bGL^{\Delta^{op}}$ such that $$\bA^n_{G_L}:=\underbrace{A\times_{G_L}\dots\times_{G_L}A}_{n+1\:factors}\ .$$
Let $\Omega^\cdot_{G_L}$ denote the de Rham complex (see \cite[3.1.2]{bss}) of the smooth stack $G_L$.  The associated chain complex of $\Omega^\cdot_{G_L}(\bA_{G_L}^\cdot)$ is a double complex with the de Rham differential $d_{dR}$ and the \v{C}ech differential $\delta$.

Note that $A\to G_L\to LX$ is an atlas. We form the simplicial object
$\bA^\cdot_{LX}\in \bLX^{\Delta^{op}}$ such that $$\bA^n_{LX}:=\underbrace{A\times_{LX}\dots\times_{LX}A}_{n+1\:factors}\ .$$
We consider the double complex $\Omega_{LX}^\cdot(\bA_{LX}^\cdot)$
Note that by \cite[Lemma 2.36]{bss} we have
$$\Omega^\cdot(X)\stackrel{Lemma \ref{twopic}}{\cong} \Gamma_X\Omega_{LX}^\cdot\cong \ker(\delta:\Omega_{LX}^\cdot(A^0_{LX})\to A^1_{LX}))\ .$$
The property that $G_L\to LX$ is a smooth gerbe with band $U(1)$ can be expressed as the fact that the diagram
$$\xymatrix{A\times_{G_L}A\ar@{=>}[r]\ar[d]&A\ar@{=}[d]\\
A\times_{LX}A\ar@{=>}[r]&A}$$
is a central $U(1)$-extension of smooth groupoids. In particular, we see that
the canonical map
$$\ker(\delta:\Omega_{LX}^\cdot(A^0_{LX})\to A^1_{LX}))\to \ker(\delta:\Omega_{LG}^\cdot(A^0_{LG})\to A^1_{LG}))
$$ is an isomorphism, i.e. we see that
\begin{equation}\label{euwndbwnw}
\Gamma_{G_L}(\Omega^\cdot_{G_L})\cong \Gamma_{LX}(\Omega^\cdot_{LX})\cong \Omega^\cdot(LX)\ .
\end{equation}

\subsubsection{}\label{ejkdekekl}

A connection on the  gerbe $f_L:G_L\rightarrow LX$ consists of a pair $(\alpha,\beta)$,
where $\alpha\in \Omega^1(A\times_{G_L}A)$ is a connection one-form on  the $U(1)$-bundle  $A\times_{G_L}A\rightarrow A\times_{LX}A$, and $\beta\in \Omega^2(A)$.
We consider $\alpha\in \Omega^1_{G_L}(\bA_{G_L}^1)$ and $\beta\in \Omega^2_{G_L}(\bA_{G_L}^0)$. 
The pair is a connection $(\alpha,\beta)$ if it satisfies:
\begin{enumerate}
\item
$\delta \beta=d_{dR}\alpha$ ,
\item
 $\delta\alpha=0$. \end{enumerate}
Note that $\delta d\beta=0$ so that there is a unique $\lambda\in \Gamma_{G_L}\Omega_{G_L}^3\stackrel{(\ref{euwndbwnw})}{\cong} \Omega^3(LX)$ which restricts to $d\beta$. We have $d\lambda=0$. 

\subsubsection{} 

Let use choose a connection $(\alpha,\beta)$, and let $\lambda\in \Omega^3(LX)$ be the associated closed three form. In \ref{ldefg} we have introduced the locally constant sheaf
$\cL$ on $LX$.  The construction \ref{dera1} gives the complex of sheaves
$$(\Omega^\cdot_{LX}(\cL)[[z]]_\lambda\ ,d_\lambda)\ .$$
Furthermore we set
$$\Omega^\cdot(LX,\cL)[[z]]_\lambda:=\Gamma_{LX}\Omega^\cdot_{LX}(\cL)[[z]]_\lambda\ .$$

\begin{ddd}\label{def36}
The delocalized $(G,\lambda)$-twisted de Rham cohomology of $X$ is defined by
$$H_{dR,deloc}^*(X,(G,\lambda)):=H^*(\Omega^\cdot(LX,\cL)[[z]]_\lambda,d_\lambda)\ .$$
\end{ddd}
In view of Lemma \ref{twopic} this is the definition  given in \cite[3.10]{math.KT/0505267}. Note that $H_{dR,deloc}^*(X,(G,\lambda))$ depends on the choice of the connection, through these groups are isomorphic for different choices (see \cite[3.11]{math.KT/0505267}).

\subsection{Comparison}

\subsubsection{}

In this subsection we prove 

\begin{theorem}\label{main1}
There is an isomorphism
$$H^*_{deloc}(X;G)\cong H^*_{dR,deloc}(X,(G,\lambda))\ .$$
\end{theorem}

Actually, this theorem follows from the following stronger statement.
Recall that $f_L:G_L\rightarrow LX$ is the pull-back of $f:G\rightarrow X$ via the canonical map $LX\to X$. Let $\R_{\bGL}\in \Sh_{\Ab}\bGL$ denote the constant sheaf with value $\R$.
\begin{theorem}\label{main2}
There is an isomorphism in $D^+(\Sh_{\Ab}\bLX)$
$$R(f_L)_*(\R_{\bG_L})\otimes_\R \cL\cong \Omega^\cdot_{LX}(\cL)[[z]]_\lambda\ .$$
\end{theorem}
The remainder of the present subsection is devoted to the proofs of Theorems \ref{main1} and \ref{main2}.

\subsubsection{}
Let us prove Theorem \ref{main2}.  First observe that
$\Omega^\cdot_{LX}(\cL)[[z]]_\lambda\cong \Omega^\cdot_{LX}[[z]]_\lambda\otimes_\R \cL$.
Therefore it suffices to show that
$$R(f_L)_*(\R_{\bGL})\cong 
 \Omega^\cdot_{LX}[[z]]_\lambda\ .$$ 
This is exactly the assertion of \cite[Theorem 1.1]{bss}, with the difference, that now $LX$ is an orbifold instead of a smooth manifold. We repeat the proof of \cite[Theorem 1.1]{bss} given by \cite[Subsection 3.2]{bss} with the following modifications (the numbers refer to the paragraphs in \cite[Subsection 3.2]{bss}:
\begin{enumerate}
\item 3.2.1 :  The manifold $X$ is replaced by the  orbifold $LX$. The gerbe $G\to X$ is replaced by the gerbe $G_L\to LX$. Furthermore, $A\to G_L$ is some atlas. It induces an atlas $A\to G_L\to LX$.
 The $U(1)$-central extension of groupoids
$(A\times_{G_L}A\Rightarrow A)\to (A\times_{LX}A\to A)$ represents a gerbe in the language of groupoids, but we can no longer refer to the paper \cite{MR1876068}.  For  existence of a connection we now refer to  \cite[Prop. 3.6]{math.KT/0505267}.
\item 3.2.2 : We use the notation $\Omega_{G_L}$ instead of $\Omega(G_L)$ for the de Rham complex of the smooth stack $G_L$.
$\Omega^\cdot(LX)$ must be interpreted as in \ref{zueuedeke}.
For the existence of connections we refer to \cite{math.KT/0505267}. The construction of the three-form associated to a connection $(\alpha,\beta)$ was explained in \ref{ejkdekekl}.
\item 3.2.6 : We must show that the map $\phi:\Omega^\cdot[[z]]_\lambda\to i^\sharp C_A(\Omega^\cdot(G_L))$ is a quasi-isomorphism. This can be shown locally. Since we can cover $LX$ by smooth manifolds the local isomorphism immediately follows from the result proved in \cite{bss}. This argument avoids repeating  the proof of \cite[Proposition 3.4]{bss}.
\end{enumerate}
\hB

\subsubsection{}

We now show Theorem \ref{main1}. We need the following well-known fact.
Let $X$ be an orbispace or orbifold and $p:X\to *$ be the projection.
Recall that $\ev\circ p_*=\Gamma_X:\Sh_{\Ab}\bX\rightarrow \Ab$. This functor
is left exact and can thus be derived. Let $\cO_X$ be the sheaf of continuous or smooth  real functions on $X$, i.e. $\cO_X=\Omega^0_X$ in the smooth case.
\begin{lem}\label{orsp}
If $F\in \Sh_{\Ab}\bX$ is a flabby sheaf and a sheaf of  $\cO_X$-modules, then
$R^i\Gamma_X (F)=0$ for $i\ge 1$.
\end{lem}
\proof
Let $A\rightarrow X$ be an orbispace (orbifold) atlas. Then $A\times_XA\Rightarrow A$ is a very proper, separated, and  {\'e}tale groupoid. Let $A^\cdot$ be the associated simplicial space (manifold). The complex $F(A^\cdot)$ represents $R\Gamma_X(F)$ by \cite[Lemma 2.41]{bss}. We now employ the method of \cite[Section 4.1]{math.GT/0508550} in order to show that $H^i(F(A^\cdot))=0$ for $i\ge 1$. We use the $\cO_X$-module structure in order to multiply by cut-off function.\hB

\subsubsection{}

We first observe that
$\Omega^\cdot_{LX}(\cL)[[z]]_\lambda$ is a complex of flabby sheaves and of $\cO_{LX}=\Omega^0_{LX}$-modules.
Therefore by Lemmas \ref{orsp} and \ref{twopic} 
we have (see  \cite[Cor. 25]{MR2172499} for a related result) $$R\Gamma_{LX} (\Omega^\cdot_{LX}(\cL)[[z]]_\lambda)  \cong \Omega^\cdot(LX,\cL)[[z]]_\lambda\ .$$
By Definition \ref{def36} the  cohomology of the right-hand side is  $H^*_{deloc,dR}(X,(G,\lambda))$.
On the other hand by Theorem \ref{main1} 
$$R\Gamma_{LX} (\Omega^\cdot_{LX}(\cL)[[z]]_\lambda)\cong
\ev \circ Rp_* \left( R(f_I)_*(\R_{\bGL})\otimes_\R \cL\right)\ .$$ Its cohomology is by \ref{umr} isomorphic to  $H^*_{deloc}(X;G)$.
\hB

\bibliographystyle{halpha.bst}

\begin{thebibliography}{BCM{\etalchar{+}}02}

\bibitem[AR03]{MR1993337}
A.~Adem and Y.~Ruan.
\newblock Twisted orbifold {$K$}-theory.
\newblock {\em Comm. Math. Phys.}, 237(3):533--556, 2003.

\bibitem[ARZ]{math.AT/0605534}
A.~Adem, Y.~Ruan, and B.~Zhang.
\newblock {A Stringy Product on Twisted Orbifold K-theory},
  arXiv:math.AT/0605534.

\bibitem[AS]{math.KT/0510674}
M.~Atiyah and G.~Segal.
\newblock {Twisted K-theory and cohomology}, arXiv:math.KT/0510674.

\bibitem[AS69]{MR0259946}
M.~F. Atiyah and G.~B. Segal.
\newblock Equivariant {$K$}-theory and completion.
\newblock {\em J. Differential Geometry}, 3:1--18, 1969.

\bibitem[AS04]{MR2172633}
M.~Atiyah and G.~Segal.
\newblock Twisted {$K$}-theory.
\newblock {\em Ukr. Mat. Visn.}, 1(3):287--330, 2004.

\bibitem[Ati61]{MR0148722}
M.~F. Atiyah.
\newblock Characters and cohomology of finite groups.
\newblock {\em Inst. Hautes \'Etudes Sci. Publ. Math.}, (9):23--64, 1961.

\bibitem[BC88]{MR928402}
P.~Baum and A.~Connes.
\newblock Chern character for discrete groups.
\newblock In {\em A f\^ete of topology}, pages 163--232. Academic Press,
  Boston, MA, 1988.

\bibitem[BCM{\etalchar{+}}02]{MR1911247}
P.~Bouwknegt, A.~L. Carey, V.~Mathai, M.~K. Murray, and D.~Stevenson.
\newblock Twisted {$K$}-theory and {$K$}-theory of bundle gerbes.
\newblock {\em Comm. Math. Phys.}, 228(1):17--45, 2002.

\bibitem[Beh04]{MR2172499}
K.~Behrend.
\newblock Cohomology of stacks.
\newblock In {\em Intersection theory and moduli}, ICTP Lect. Notes, XIX, pages
  249--294 (electronic). Abdus Salam Int. Cent. Theoret. Phys., Trieste, 2004.

\bibitem[BS]{math.GT/0508550}
U.~Bunke and T.~Schick.
\newblock {T-duality for non-free circle actions}, arXiv:math.GT/0508550.

\bibitem[BSSa]{bssf}
U.~Bunke, Th. Schick, and M.~Spitzweck.
\newblock {Foundations of sheaf theory on topological stacks}.
\newblock In preparation.

\bibitem[BSSb]{bssm}
U.~Bunke, Th. Schick, and M.~Spitzweck.
\newblock {$T$-duality and periodic twisted cohomology}.
\newblock In preparation.

\bibitem[BSSc]{bss}
U.~Bunke, T.~Schick, and M.~Spitzweck.
\newblock {Sheaf theory for stacks in manifolds and twisted cohomology for
  $S^1$-gerbes}, arXiv:math.KT/0603698.

\bibitem[BX]{math.DG/0605694}
K.~Behrend and P.~Xu.
\newblock {Differentiable Stacks and Gerbes}, arXiv:math.DG/0605694.

\bibitem[CR04]{MR2104605}
W.~Chen and Y.~Ruan.
\newblock A new cohomology theory of orbifold.
\newblock {\em Comm. Math. Phys.}, 248(1):1--31, 2004.

\bibitem[Cra]{math.DG/0008064}
M.~Crainic.
\newblock {Differentiable and algebroid cohomology, van Est isomorphisms, and
  characteristic classes}, arXiv:math.DG/0008064.

\bibitem[FHT]{math.AT/0312155}
D.~S. Freed, M.~J. Hopkins, and C.~Teleman.
\newblock Twisted $k$-theory and loop group representations i.

\bibitem[Hei05]{heinloth}
J.~Heinloth.
\newblock Survey on topological and smooth stacks.
\newblock In {\em Mathematisches Institut G{\"o}ttingen, WS04-05 (Y. Tschinkel,
  ed.)}, pages 1--31. 2005.

\bibitem[Hit01]{MR1876068}
N.~Hitchin.
\newblock Lectures on special {L}agrangian submanifolds.
\newblock In {\em Winter School on Mirror Symmetry, Vector Bundles and
  Lagrangian Submanifolds (Cambridge, MA, 1999)}, volume~23 of {\em AMS/IP
  Stud. Adv. Math.}, pages 151--182. Amer. Math. Soc., Providence, RI, 2001.

\bibitem[Hov99]{MR1650134}
M.~Hovey.
\newblock {\em Model categories}, volume~63 of {\em Mathematical Surveys and
  Monographs}.
\newblock American Mathematical Society, Providence, RI, 1999.

\bibitem[Kaw78]{MR0474432}
T.~Kawasaki.
\newblock The signature theorem for {$V$}-manifolds.
\newblock {\em Topology}, 17(1):75--83, 1978.

\bibitem[Kaw81]{MR641150}
T.~Kawasaki.
\newblock The index of elliptic operators over {$V$}-manifolds.
\newblock {\em Nagoya Math. J.}, 84:135--157, 1981.

\bibitem[LU]{math.AT/0307114}
E.~Lupercio and B.~Uribe.
\newblock {Holonomy for Gerbes over Orbifolds}, arXiv:math.AT/0307114.

\bibitem[LU02]{MR1950946}
E.~Lupercio and B.~Uribe.
\newblock Loop groupoids, gerbes, and twisted sectors on orbifolds.
\newblock In {\em Orbifolds in mathematics and physics (Madison, WI, 2001)},
  volume 310 of {\em Contemp. Math.}, pages 163--184. Amer. Math. Soc.,
  Providence, RI, 2002.

\bibitem[LU04]{MR2045679}
E.~Lupercio and B.~Uribe.
\newblock Gerbes over orbifolds and twisted {$K$}-theory.
\newblock {\em Comm. Math. Phys.}, 245(3):449--489, 2004.

\bibitem[LUX]{math.AT/0512658}
E.~Lupercio, B.~Uribe, and M.~A. Xicotencatl.
\newblock {Orbifold String Topology}, arXiv:math.AT/0512658.

\bibitem[Met]{math.DG/0306176}
D.~Metzler.
\newblock {Topological and Smooth Stacks}, arXiv:math.DG/0306176.

\bibitem[MM97]{Minasian:1997mm}
R.~Minasian and G.~W. Moore.
\newblock K-theory and ramond-ramond charge.
\newblock {\em JHEP}, 11:002, 1997, hep-th/9710230.

\bibitem[MS03]{MR1977885}
V.~Mathai and D.~ Stevenson.
\newblock Chern character in twisted {$K$}-theory: equivariant and holomorphic
  cases.
\newblock {\em Comm. Math. Phys.}, 236(1):161--186, 2003.

\bibitem[Noo]{math.AG/0503247}
B.~ Noohi.
\newblock {Foundations of Topological Stacks I}, arXiv:math.AG/0503247.

\bibitem[Rua]{math.AG/0005299}
Y.~Ruan.
\newblock {Discrete torsion and twisted orbifold cohomology},
  arXiv:math.AG/0005299.

\bibitem[Rua02]{MR1941583}
Y.~Ruan.
\newblock Stringy geometry and topology of orbifolds.
\newblock In {\em Symposium in Honor of C. H. Clemens (Salt Lake City, UT,
  2000)}, volume 312 of {\em Contemp. Math.}, pages 187--233. Amer. Math. Soc.,
  Providence, RI, 2002.

\bibitem[Tu99]{MR1671260}
J.~L. Tu.
\newblock La conjecture de {N}ovikov pour les feuilletages hyperboliques.
\newblock {\em $K$-Theory}, 16(2):129--184, 1999.

\bibitem[TXa]{math.KT/0505267}
J.~L. Tu and P.~Xu.
\newblock {Chern character for twisted K-theory of orbifolds},
  arXiv:math.KT/0505267.

\bibitem[TXb]{math.KT/0604160}
J.~L. Tu and P.~Xu.
\newblock {The ring structure for equivariant twisted K-theory},
  arXiv:math.KT/0604160.

\bibitem[TXLG04]{MR2119241}
J.~L.Tu, P.~Xu, and C.~Laurent-Gengoux.
\newblock Twisted {$K$}-theory of differentiable stacks.
\newblock {\em Ann. Sci. \'Ecole Norm. Sup. (4)}, 37(6):841--910, 2004.

\bibitem[Wit98]{Witten:1998cd}
E.~Witten.
\newblock D-branes and K-theory.
\newblock {\em JHEP}, 12:019, 1998, hep-th/9810188.

\end{thebibliography}

\newcommand{\etalchar}[1]{$^{#1}$}

\end{document}